\documentclass[11pt]{article}
\usepackage{amssymb}
\usepackage{amsmath, amscd}

\usepackage[marginratio=1:1]{geometry}
\usepackage{layout}
\usepackage{dsfont}
\usepackage{enumerate}
\usepackage{graphicx}
\usepackage{tikz-cd}

\usepackage{tikz}
\usetikzlibrary{calc, intersections}

\newtheorem{lem}{Lemma} 
\newtheorem{rem}{Remark}

\newtheorem{thm}{Theorem}[section]
\newtheorem{prop}[thm]{Proposition}

\newtheorem{cor}[thm]{Corollary}
\newtheorem{dfn}[thm]{Definition}

\newcommand{\pf}{\noindent{\bf Proof.}\ }
\newcommand{\qed}{\begin{flushright}$\Box$\ \ \ \ \ \
	\end{flushright}}

\usepackage{tikz}
\usetikzlibrary{calc, intersections}

\title{Classifying Linear Canonical Relations}
\author{Jonathan Lorand, ETH Zurich \\
Summer 2014\\
\\
Master's Thesis in Mathematics under the advisorship of \\
Giovanni Felder, ETH Zurich, and Alan Weinstein, UC Berkeley\\
}
\date{}

\begin{document}
\maketitle

\

\begin{abstract}
We consider the problem of classifying, up to conjugation by linear symplectomorphisms, linear canonical relations (lagrangian correspondences) from a finite-dimensional symplectic vector space to itself. We give an elementary introduction to the theory of linear canonical relations and present partial results toward the classification problem. This exposition should be accessible to undergraduate students with a basic familiarity with linear algebra. 
\end{abstract}
\

\newpage

\tableofcontents

\newpage

\section{Introduction}

\subsection{Summary and Context}

The subject of the present work is the problem of classifying, up to linear symplectomorphism, linear canonical relations from a symplectic vector space $(V,\omega)$ to itself. We begin with a precise formulation of this classification problem and a short overview of this work's contents, which comprise partial results toward the solution of this classification problem, as well as an exposition of the tools and results used toward this aim. 

The basic setting involves a finite-dimensional real vector space $V$, equipped with a symplectic form $\omega$, \text{i.e.} a non-degenerate antisymmetric bilinear form $\omega: V \times V \rightarrow \mathbb{R}$. If one multiplies $\omega$ by $-1$ this also gives a symplectic form; we use the notation $V^{-}$ to indicate when $V$ is endowed with this symplectic structure. On the space $V \oplus V^{-}$ we have a naturally defined symplectic form which is the direct sum of the symplectic forms on $V$ and $V^{-},$ 
$$
((v,\bar{v}),(u,\bar{u})) \longmapsto \omega(v,u) - \omega(\bar{v},\bar{u}) \quad \quad  v,u \in V, \bar{v}, \bar{u} \in V^-
$$
A linear subspace $L \subset V\oplus V^{-}$ is a \textbf{linear canonical relation} if it is a \textbf{lagrangian} subspace, i.e. such that the symplectic form on $V\oplus V^{-}$ is identically zero when restricted to $L$ and such that $L$ has the maximum possible dimension for such a subspace, that is $\dim L = \dim V$. If $(V,\omega)$ and $(\hat{V},\hat{\omega})$ are symplectic vector spaces and $L$ and $\hat{L}$ linear canonical relations in $V \oplus V^-$ and $\hat V \oplus \hat V^-$ respectively, we say they are \textbf{equivalent} if there exists a linear symplectomorphism $S: V \rightarrow \hat V$ (a linear isomorphism satisfying $\hat \omega(Sv,Su)) = \omega(v,u)$) such that 
$$
(v,u) \in L  \ \Leftrightarrow \ (Sv,Su) \in \hat{L}
$$
The classification problem at hand is to determine invariants which uniquely determine the equivalence classes of this equivalence relation and to give normal forms for these, \text{i.e.} unique representatives of each equivalence class. 

The graph of any linear symplectomorphism from $V$ to itself is a special case of a linear canonical relation. If we consider only such linear canonical relations, the above classification problem amounts to the problem of giving normal forms for the conjugacy classes of the group of linear symplectomorphisms on $\mathbb{R}^{2n}$. This problem has various solutions and a long history, which extends from Williamson \cite{Williamson} to, most recently, Gutt \cite{Gutt}. The paper of Laub and Meyer \cite{Laub Meyer} contains a helpful albeit somewhat dated survey of this history. In \cite{Towber} Towber carries out the classification of linear relations up to linear isomorphism, \text{i.e.} the more general version of our present classification problem which does not include any symplectic structure. We will discuss some basic properties of linear relations and Towber's results in Section \ref{Linear Relations}. In Section \ref{Linear Canonical Relations} we describe basic properties of linear canonical relations, the most important being the result of Benenti and Tulczyjew stating that a linear canonical relation is described by two coisotropic subspaces and an induced linear symplectomorphism between their reduced spaces. In particular, as a first step in our classification problem it is necessary to classify pairs of coisotropic subspaces up to linear symplectomorphism. This is done in Section \ref{Coisotropic Pairs}. Finally in Section \ref{Classification of Linear Canonical Relations} we present partial results towards a full solution of the classification problem. 

Linear canonical relations are the linear counterparts of canonical relations (also known as lagrangian correspondences), which are lagrangian submanifolds of the product of two symplectic manifolds. One motivation for the study of such objects arises in the context of quantization, where one is interested in making rigorous the correspondences between mathematics associated to classical mechanics on the one hand, and to quantum mechanics on the other. On the classical side a symplectic manifold, symplectomorphisms of and functions on this manifold are used to describe the possible states, the symmetries, time evolution, and the observables of a physical system. The corresponding objects on the quantum side are usually given by a Hilbert space together with an algebra of operators on this space and probability measures constructed from such operators. One approach to the quantization problem is to try to find appropriate descriptions on both the classical and quantum sides which exhibit the structure of a category, and for which the passage from classical to quantum is functorial. 

In this context, canonical relations are a way to generalize the idea of a symplectomorphism with the goal of obtaining a category in which they are the morphisms. In the linear case, the linear canonical relations do in fact define the morphisms of a category of which the objects are symplectic vector spaces. In the non-linear case, where the objects are symplectic manifolds, this is not so. Both the linear and the non-linear cases exhibit delicacies when one attempts to compose canonical relations which do not satisfy a certain transversality condition. In the linear situation this manifests itself in that the composition operation becomes non-continuous in terms of the usual topology on the space of all lagrangian subspaces. A solution to this problem is presented by Li-Bland and Weinstein in  \cite{Li-Bland Weinstein}. In the non-linear case, further difficulties arise. Although these issues are beyond the scope of the present exposition, they form part of a greater context in which it is embedded. For more on the these general topics, we refer the interested reader to \cite{G of Q} and \cite{Symplectic Categories}, and the many references therein. 

\subsection{Conventions and Notation}

Throughout the text, notation is usually standard or is introduced as needed. We note here though that everything takes place in finite dimensions and that everything is linear.  Many times the adjective ``linear'' will thus be omitted (and tacitly implied). In this vein, the terms ``linear symplectomorphism", ``symplectomorphism" and ``symplectic map" will be used as synonyms. 

\subsection{Acknowledgments}

I wish to thank first and foremost Alan Weinstein, for his generous kindness, humor, patience and guidance throughout my stay at UC Berkeley, of which this thesis is but one product. My gratitude also goes to Giovanni Felder, for his encouragement, mentorship and the reading course last fall, which led to where I am now. Lastly, I express my appreciation toward the Anna and Hans K\"agi Foundation, the ETH Zurch, and the City of Zurich for their financial sponsorship, and to my family and friends for their support. 

\newpage

\section{Symplectic Linear Algebra}\label{Symplectic Linear Algebra}

\subsection{Basic notions}\label{Basic notions}

We recount here some linear algebraic definitions and constructions in the setting of a symplectic vector space $(V,\omega)$. If $E,F \subset V$ are subspaces, we say they are $\omega$-\textbf{orthogonal} if $ \omega(e,f) = 0$ for all $e \in E, f \in F$. For any subspace $W \subset V$, \emph{its} $\omega$-orthogonal subspace is 
$$
W^{\omega}:= \{ v \in V \mid \omega(v,u) = 0 \ \forall u \in W \}
$$
This space is not in general a complement of $W$ in $V$, but its dimension is complementary
$$
\dim W + \dim W^{\omega} = \dim V.
$$
One way to see this is via the mapping
$$
\tilde{\omega} : V \rightarrow V^* \quad v \mapsto \omega(v, \cdot).
$$
The non-degeneracy of $\omega$ means that $\tilde \omega$ is an isomorphism. By post-composing this map with the restriction $V^* \rightarrow U^*$ one obtains an epimorphism $V \rightarrow U^*$ with kernel $U^{\omega}$, and consequently an isomorphism $V/U^{\omega} \simeq U^*$.

As an operation on subspaces, taking the $\omega$-orthogonal is involutive and exchanges sums and intersections
$$
(W^{\omega})^{\omega} = W \quad \quad (E \cap F)^{\omega} = E^{\omega} + F^{\omega} \quad \quad \quad W,E,F \subset V.
$$
If we restrict $\omega$ to $W$, its kernel is $W \cap W^{\omega}$. One defines 
$$
\text{rank} (W) := \dim W - \dim W \cap W^{\omega}.
$$
A subspace $W \subset V$ is called \textbf{symplectic} if the restriction $\omega_W$ of $\omega$ to $W$ defines a symplectic form, making $(W,\omega_W)$ a symplectic space in its own right. That $\omega_W$ be non-degenerate means $\text{rank} (W) = \dim W$ and $W$ is symplectic if and only if $W \cap W^{\omega} = 0$. In particular, if $W$ is symplectic then $W^{\omega}$ is too, and one has an $\omega$-orthogonal decomposition $V = W \oplus W^{\omega}$. Further fundamental types of subspaces can also be defined via the relation of a space with its orthogonal: $W$ is 
\begin{align*}
\textbf{istropic} & \quad \text{if } W \subset W^{\omega} \\
\textbf{coistropic} & \quad \text{if } W^{\omega} \subset W \\
\textbf{lagrangian} & \quad \text{if } W^{\omega} = W  
\end{align*}
In particular, $W$ is lagrangian if and only if it is isotropic (or coisitropic) and $\dim W = \dim W^{\omega}$. 

A \textbf{linear symplectomorphism} or \textbf{symplectic map} from one symplectic space $(V,\omega)$ to another $(\hat V,\hat \omega)$ is a linear isomorphism $S:V \rightarrow \hat V$ such that $\hat \omega (Sv,Su) = \omega(u,v)$. For such a map one has
\begin{align*}
S(W)^{\hat \omega} & = \{\hat v \in \hat V \mid  \hat \omega (\hat v,Su) = 0 \ \forall u \in W\} \\
				& = \{\hat v \in \hat V \mid \omega (S^{-1}\hat v,u) = 0 \ \forall u \in W\} \\
				& = \{Sv \in \hat V \mid \omega (v,u) = 0 \ \forall u \in W\}  \\
				& = S(W^{\omega})
\end{align*}
hence in particular $S(W)$ is symplectic/(co)isotropic/lagrangian when $W$ is, and any $\omega$-orthogonal decomposition $V = E \oplus F$ is mapped to an $\hat \omega$-orthogonal decomposition $S(E)\oplus S(F) = \hat V$. In general, if $V = E \oplus F$ and $\hat V = \hat E \oplus \hat F$ are any given decompositions (possibly not  $\omega$-orthogonal), we say a map $S:V \rightarrow \hat V$ satisfying $S(E) = \hat E$ and $S(F) = \hat F$ \textbf{respects the decompositions} in $V$ and $\hat V$. If $E$, $F$, $\hat E$ and $\hat F$ are symplectic subspaces and $S$ is a symplectic map which respects the decompositions, then $S\vert_E: E \rightarrow \hat E$ and $S \vert_F : F \rightarrow \hat F$ are again symplectic maps. The following gives a kind of converse to this fact. 

\begin{lem} \label{symplectic direct sum} Let $V = E \oplus F$  and $V = \hat E \oplus \hat F$ be two $\omega$-orthogonal direct sum decompositions with symplectic subspaces. Let $\phi: E \rightarrow \hat E$ and $\psi : F \rightarrow \hat F$, and set $\sigma = \phi \oplus \psi : V \rightarrow \hat V$. Then $\sigma$ is symplectic iff $\phi$ and $\psi$ are symplectic.
\end{lem}

\pf
Assume $\phi$ and $\psi$ are symplectic, and let $w = w_1 + w_2$ denote the decomposition of any $w \in V$ with respect to the splitting $E \oplus F$. For $u,v \in V$ one has 
\begin{align*}
\hat \omega(\sigma u, \sigma v) & = \hat \omega(\phi u_1 + \psi u_2 , \phi v_1 + \phi v_2) \\
						& = \hat \omega(\phi u_1, \phi v_1 ) + \hat \omega(\psi u_2 , \psi v_2) \\
						& = \omega(u_1,v_1) + \omega(u_2,v_2) \\
						& = \omega(u_1 + u_2, v_1 + v_2) \\
						& = \omega(u,v)
\end{align*}
where the second and fourth equalities hold due to the orthogonality of $\hat E$ and $\hat F$, and $E$ and $F$, respectively. This shows that $\sigma$ is symplectic. The converse statement, when $\sigma$ is assumed symplectic, is clear since then $\phi = \sigma \vert_E$ and $\psi = \sigma \vert_F$. 
 
\qed
In the above, the orthogonality condition on $E$ and $F$ (and $\hat E$ and $\hat F$) amounts to $E \oplus F$ being naturally symplectomorphic to the external direct sum of two separate symplectic spaces $(E,\omega_E)$ and $(F, \omega_F)$, endowed with the direct sum symplectic form $\omega_E \oplus \omega_F$ defined by $((e,f),(e',f')) \mapsto \omega_E(e,e') + \omega_F(f,f')$. 

The 2-dimensional symplectic spaces are the basic building blocks of symplectic vector spaces in the sense that any such may be decomposed as 
$$
V = \bigoplus_{i=1}^{n} E_i  
$$
where each $E_i$ is a 2-dimensional symplectic subspace and $n \in \mathbb{N}$. In particular $\dim V = 2n$ is even. We will henceforth always use $n$ to denote half the dimension of $V$. To prove the above, one can use an iterative Gram-Schmidt-type process to construct a basis $\{q_1,...,q_n, p_1,...,p_n \}$ of $V$ such that 
\begin{equation}\label{symplectic basis condition}
\omega(q_i,p_j) = \delta_{ij}, \quad \omega(q_i,q_j) = \omega(p_i,p_j) = 0  \quad  \quad \forall i,j \in \{1,....,n\} \quad \quad \quad 
\end{equation}
Any ordered basis satisfying this property is a \textbf{symplectic basis}. For such a basis, the 2-dimensional subspaces $E_i = \langle q_i, p_i \rangle $ are $\omega$-orthogonal when $i \neq j$, and the matrix associated to each $\omega_{E_i}$ 
$$
\left(  
\begin{array}{cc}
\omega(q_i,q_i) & \omega(q_i,p_i)  \\
\omega(p_i,q_i) & \omega(p_i,p_i)
\end{array}
\right) 
=
\left(  
\begin{array}{cc}
0 & 1 \\
-1& 0
\end{array}
\right) 
$$
is non-singular, so each $E_i$ is a symplectic subspace. To construct a symplectic basis, begin with any vector $q_1 \in V$. Since $\omega$ is non-degenerate there exists a vector $p_1$ such that $\omega(q_1,p_1) \neq 0$, and this vector may be normalized if need be so that $\omega(q_1,p_1) = 1$. Note that $\omega(q_1,p_1) \neq 0$ implies that $q_1$ and $p_1$ are linearly independent since by antisymmetry $\omega(v,v) = 0$ for any $v \in V$. The subspace $E_1 = \langle q_1, p_1 \rangle$ is a symplectic subspace, and its $\omega$-orthogonal complement is such that $V = E_1 \oplus E_1^{\omega}$. If this complement is zero we are done, otherwise we can choose any vector $q_2 \in E_1^{\omega}$ and, because $\omega_{E_1^{\omega}}$ is non-degenerate, we can also find $p_2$ such that $\omega(q_2,p_2)= 1$. In particular the dimension of $E_1^{\omega}$, if non-zero, must be at least two. One sets $E_2 := \langle q_2,p_2 \rangle $ and continues this process in the $\omega$-orthogonal complement of the symplectic space $E_1 \oplus E_2$. Proceeding iteratively one must reach a point after a finite number $n \in \mathbb{N}$ iterations where the $\omega$-orthogonal complement of $ E_1 \oplus ... \oplus E_n$ is zero, since the dimension of such complements decreases by 2 with each step, and each complement has dimension either greater than 1 or equal to 0. Hence the $q_i$ and $p_i$ span $V$, and they are linearly independent and satisfy (\ref{symplectic basis condition}) by construction. The fact that every symplectic vector space admits a symplectic basis is the linear version of Darboux's theorem.

\subsection{Lagrangian splittings}\label{Lagrangian splittings}

In addition to direct sum decompositions into symplectic subspaces, one may also consider decompositions of the form $V=L_1 \oplus L_2$, where $L_1$ and $L_2$ are lagrangian subspaces. In such a case we call $(L_1,L_2)$ a \textbf{lagrangian splitting} of $V$ or a \textbf{transverse lagrangian pair}. Because a symplectic basis always exists, so also do lagrangians and lagrangian splittings. If $\{q_1,...,q_n,p_1,...p_n\}$ is a symplectic basis, $\langle q_1,...,q_n \rangle$ and $\langle p_1,...,p_n \rangle $ are lagrangian subspaces forming a lagrangian splitting. 

\begin{prop} \label{lagrangian complements}Let $L_1 \subset V$ be lagrangian. Then there exists a lagrangian subspace $L_2 \subset V$ such that $(L_1,L_2)$ is a lagrangian splitting of $V$. 
\end{prop}

\pf 
Let $\dim V = 2n$ and let $\{q_1,...,q_n\}$ be a basis of $L:=L_1$. We proceed in a similar manner as in the iterative construction of a symplectic basis, though now the $q_i$ must be in the prescribed subspace $L$. To $q_1$ we can find $p_1 \notin L$ such that $\omega(q_i,p_1) = \delta_{i1}, \ i \in \{1,...,n\},$ by choosing $p_1$ in $\langle q_2,...,q_n \rangle ^{\omega} \backslash \langle q_1 \rangle ^{\omega}$ and scaling if necessary. This is possible because $\langle q_2,...,q_n \rangle ^{\omega}$ is $n+1$-dimensional, and $ \langle q_2,...,q_n \rangle ^{\omega}  \cap \langle q_1 \rangle ^{\omega} = L^{\omega} = L$ is $n$-dimensional; $p_1 \notin L$ is clear because $L \subset \langle q_1 \rangle ^{\omega}$. Set $E_1 = \langle q_1,p_1 \rangle$. Because $q_1,p_1 \in \langle q_2,...,q_n \rangle ^{\omega}$, we have $\langle q_2,...,q_n \rangle \subset E_1^{\omega} \cap L$, and this is an equality for dimension reasons, since $q_1 \in E_1 \cap L$. So $L_2 := \langle q_2,...,q_n \rangle$ is a Langrangian subspace of the symplectic space $E_1^{\omega}$ and one can iterate the procedure by choosing $p_2 \in ( \langle q_3,...,q_n \rangle ^{\omega} \cap E_1^{\omega} )\backslash (\langle q_2 \rangle ^{\omega} \cap E_1^{\omega} ) $ and normalizing such that $\omega(q_2,p_2) = 1$. After a finite number of steps one will have constructed a symplectic basis $\{q_1,...,q_n,p_1,...p_n\}$, which means that $L_2:= \langle p_1,...,p_n \rangle$ is a langrangian complement of $\langle q_1,...,q_n \rangle = L$.
\qed

In the proof above, we constructed a symplectic basis which extends a given basis of $L_1$. One can in prescribe a lagrangian $L_1$, a basis   of $L_1$, as well as a lagrangian complement $L_2$. 

\begin{prop} \label{lagrangian pair dual basis} Let $(L_1, L_2)$ be a lagrangian splitting of V, and let $\{q_1,...,q_n \}$ be a basis of $L_1$. Then there exists exactly one basis $\{p_1,...,p_n\}$ of $L_2$ such that $\{q_1,...,q_n,p_1,...,p_n\}$ is a symplectic basis of $V$. 
\end{prop}

\pf One way to show this is to note that the isomorphism $\tilde \omega : V \rightarrow V^*$ induces an isomorphism $L_1 \rightarrow L_2^*$, and the images of the $q_i$ basis vectors in $L_2$ give the dual basis of the desired basis $\{p_i \}$ of $L_2$. Indeed, by precomposing $\tilde \omega$ with the inclusion map $L_1 \hookrightarrow V$ and post-composing with the restriction $V^* \rightarrow L_2^*, \ l \mapsto l\vert_{L_2}$, one obtains a linear mapping $\tilde \omega_{12} : L_1 \rightarrow L_2^*$. We check that it is an isomorphism. 

For surjectivity, consider any $l_2 \in L_2^*$. We can extend $l_2$ to a linear map $l \in V^*$ by setting $l = l_2 \circ \pi_2$, where $\pi_2 : V \rightarrow L_2$ is the projection onto the second lagrangian subspace with respect to the splitting $V = L_1 \oplus L_2$. Then $l$ has a unique pre-image $v = \tilde \omega^{-1} (l) \in V$ such that $l(w) = \omega (v,w) \ \forall \ w \in V$. But $v$ is in fact in $L_1$: If $w=w_1 + w_2$ is the decomposition of an arbitrary $w \in V$ into its $L_1$ and $L_2$ components (and $v_2$ the $L_2$ component of $v$), then we find $\omega (v_2,w) = \omega(v_2,w_1) + \omega(v_2,w_2) = 0$ (the first term is $l(w_1) = l_2(\pi_2(w_1)) = 0$ and the second term vanishes because $L_2$ is lagrangian). Hence $v_2=0$, i.e. $v\in L_1$, and $\tilde \omega_{12}$ maps $v$ to $l_2$. 

For injectivity, consider $v \in L_1$ in the kernel of $\omega_{12}$. Then for arbitrary $w\in V$, $\omega (v,w) = \omega(v,w_1) + \omega(v,w_2) = 0$: the first term vanishes by the lagrangian property and the second term because  $ \omega(v,w_2) = \tilde \omega_{12}(v)(w) = 0$ by the assumption $v\in \ker{ \tilde \omega_{12}}$. Since $w$ was arbitrary, $v = 0$. 

So $\tilde \omega_{12}$ maps the basis $\{q_1,...,q_n \}$ to a basis of $L_2^*$. Let $\{ p_1,...,p_n\}$ be the dual basis in $L_2$ of this basis, i.e. such that 
$$\tilde \omega_{12}(q_i)(p_j)= \delta_{ij} \ \forall \  1 \leq i,j \leq n \quad (*)$$ 
Because $\tilde \omega_{12}(q_i)(p_j) = \omega(q_i, p_j)$, this condition - together with the fact that the $q_i$ and $p_j$ each span a lagrangian subspace - means precisely that $\{q_1,...,q_n,p_1,...,p_n\}$ is a symplectic basis of $V$. It is clear that the subbasis $\{p_1,...,p_n\}$ is unique given $\{q_1,...,q_n \}$ and $L_2$, since any such has to fulfill $(*)$, which has a unique solution. 
\qed

\begin{prop} \label{symplectomorphic lagrangian pairs}
Let $(L_1,L_2)$ be a lagrangian splitting of $V$. Set $T = \mathbb{R}^n$ and equip $T\times T^*$ with the symplectic form $\overline \Omega((v,\alpha),(w,\beta)) =  \beta(v) - \alpha (w)$. Then there exists a symplectomorphism $\phi: T \times T^* \rightarrow V$ such that $\phi (T) = L_1$ and $\phi (T^*) = L_2$.
\end{prop}
 
\pf
A canonical symplectic basis for $T \times T^*$ is given by the standard basis on $\mathbb{R}^n$, together with its dual basis. By Proposition \ref{lagrangian pair dual basis}, we can find a symplectic basis $\{ q_1,...,q_n,p_1,...,p_n \}$ in $V$ such that $\langle q_1,...,q_n \rangle = L_1$ and $\langle p_1,...,p_n \rangle = L_2$. Then the mapping defined by $e_i \mapsto q_i$ and $p_i \mapsto e_i^*$ gives a symplectomorphism $\phi$ as desired. 
\qed

\begin{cor}
Let $V$ and $\hat V$ be two symplectic vector spaces of the same dimension, and $V = L_1 \oplus L_2$ and $\hat V = \hat L_1 \oplus \hat L_2$ decompositions into complementary lagrangian pairs. Then there exists a symplectomorphism $S:V \rightarrow \hat V$ such that $S(L_1) = \hat L_1$ and $S(L_2) = \hat L_2 $.
\end{cor}

\subsection{Reduction, Witt-Artin decomposition}\label{Reduction, Witt-Artin decomposition}

The quotient construction known as \textbf{symplectic reduction} produces a symplectic space from any subspace $W \subset V$. 

\begin{lem}
Let $W \subset V$ be any subspace. Then $ W / (W \cap W^{\omega})$ carries a natural induced symplectic structure $[\omega]$, given by $[\omega] ([u],[v]) := \omega(u,v)$ for all $u,v \in W$. 
\end{lem}

\pf
To check that the form $[\omega]$ is well-defined, let $u,\tilde{u},v,\tilde{v} \in W$ be such that $[u] = [\tilde{u}]$ and $[v]=[\tilde{v}]$ . Then $u = \tilde{u} + k$ and $v = \tilde{v} + l$ for some $k,l \in W \cap W^{\omega}$ and 
$$
\omega(u,v) = \omega(\tilde{u},v) + \omega(\tilde{u},l) + \omega(k,\tilde{v}) + \omega(k,l)
$$
is equal to $\omega(\tilde{u},\tilde{v})$ because the three right-hand terms above vanish since $k$ and $l$ are $\omega$-orthogonal to all of $W$. To see that $[\omega]$ is non-degenerate, assume $[u]$ is such that $[\omega]([u],[v]) = 0$ for all $[v] \in W / (W \cap W^{\omega})$. This implies $\omega(u,v) = 0$ for all $v \in W$, so $u \in W^{\omega} \cap W$ and hence $[u] = [0]$. 
\qed

The reduced space $ W / (W \cap W^{\omega})$ will sometimes be denoted $V^{W}$ and $\rho : W \rightarrow V^{W}, w \mapsto [w]$ is the reduction map associated to the reduction of $V$ by $W$. In the special case when $W$ is a coisotropic subspace, this map has the following useful property. 

\begin{lem}\label{coisotropic reduction}
Let $W \subset V$ be a coisotropic subspace. If $L \subset V$ is an isotropic (coisotropic) subspace, then the image of $L \cap W$ under $\rho$ is an isotropic (coisotropic) subspace of $V^{W}$. In particular, when $L$ is lagrangian, then so is $\rho(L\cap W)$.
\end{lem}

\pf
Let $[U]$ denote the image under $\rho$ of any subspace $U \subset W$. One has 
\begin{align*}
[U \cap W]^{[\omega]} & = \{ [v] \in V^{W} \mid \omega_W(v,l) = 0 \ \forall l \in U \cap W \} \\
				& = [(U \cap W)^{\omega_W}] = [(U \cap W)^{\omega} \cap W] \\
				& = [(U^{\omega} + W^{\omega}) \cap W] \\
				& = [U^{\omega} \cap W]  	
\end{align*}
where the last equality holds because $W^{\omega} \subset W$. Since the partial ordering of inclusion is preserved under the map $\rho$ (signified by brackets), we see that $[U^{\omega} \cap W]  \subset [U \cap W]^{[\omega]}$ if $U$ is isotropic, and $[U^{\omega} \cap W]  \supset [U \cap W]^{[\omega]}$ is $U$ is coisotropic. 
\qed

For any subspace $W \subset V$, the subspace $W \cap W^{\omega}$ is the kernel of $\omega_{W^{\omega}}$ as well as $\omega_W$, hence one can in a sense simultaneously perform a reduction with respect to both $W$ and $W^{\omega}$. Lifting back to $V$, this induces a decomposition of $V$ as an $\omega$-orthogonal direct sum of symplectic subspaces. 

\begin{prop}[Witt-Artin decomposition] \label{Witt-Artin decomposition}
 Let $W \subset V$ be any subspace, and $E$ and $F$ complements of $W\cap W^{\omega}$ in $W$ and $ W^{\omega}$ respectively. Then $E$ and $F$ are symplectic subspaces and $\omega$-orthogonal, and $V$ decomposes as the $\omega$-orthogonal direct sum 
$$ V = E \overset{\omega}{\oplus} F \overset{\omega}{\oplus} (E \oplus F)^{\omega}$$ 
Moreover, $W\cap W^{\omega}$ is a lagrangian subspace of $(E \oplus F)^{\omega}$. 
\end{prop}

\pf
Let $\pi: W \rightarrow E$ be the projection map associated to the decomposition $W = W\cap W^{\omega} \oplus E$. This induces an isomorphism $\tilde \pi : W/W\cap W^{\omega} \rightarrow E$ such that $\tilde \pi ([v]) = \pi(v)$ for all $v\in W$. Under this map, the symplectic form $[\omega]$ on the reduced space $W/W\cap W^{\omega}$ is pushed forward to a symplectic form on $E$, and $\tilde \pi_*[\omega] = \omega_E$:  
$$\tilde \pi_*[\omega](e_1,e_2) = [\omega](\tilde \pi_E^{-1}e_1,\tilde \pi^{-1}e_2) = [\omega]([e_1],[e_2]) = \omega(e_1,e_2) \quad \quad \forall e_1,e_2 \in E $$
Thus $E$ is symplectic, and by analogous arguments $F$ is symplectic as well. $E$ and $F$ are $\omega$-orthogonal because $E \subset W$ and $F \subset W^{\omega}$. As a consequence, $E \cap F = 0$ and $E \oplus F$ is symplectic also. From this it follows that $V = E \oplus F \oplus (E \oplus F)^{\omega}$. 

Finally, $W\cap W^{\omega}$ is in $(E \oplus F)^{\omega}$ since it is in $E^{\omega}$ and $F^{\omega}$ each, and $(E + F)^{\omega} = E^{\omega} \cap F^{\omega}$. Clearly $W\cap W^{\omega}$ is isotropic. To see that it is lagrangian in $G:=(E + F)^{\omega}$, note that $W^{\omega} + W = E \oplus F \oplus (W \cap W^{\omega})$ and recall the general fact that if $U,X,Y$ are subspaces such that $U \supset X$ and $X \cap Y = 0$,  it holds that $U \cap (X \oplus Y) = X \oplus (U \cap Y)$. We now calculate 
\begin{align*}
(W\cap W^{\omega})^{\omega \vert_G} & = (W\cap W^{\omega})^{\omega} \cap G \\ 
& = [(E \oplus F) \oplus (W \cap W^{\omega})] \cap G \\ 
& = [(E \oplus F) \cap G] \oplus [W \cap W^{\omega}] \\ 
& = (W\cap W^{\omega})
\end{align*}
where the last inequality uses the fact that $(E\oplus F) \cap (E \oplus F)^{\omega} = 0$ and the second to last uses the general fact about subspace above, with $G$ in the role of $U$. 
\qed

\begin{center}
\begin{tikzpicture}[scale = 1]

\path[name path=V, draw] (2.5,2.5) circle [radius=2];

\path[name path=V, clip] (2.5,2.5) circle [radius=2];

\coordinate (c) at (2.5,2.5);
\coordinate (r1) at (0,2.5);
\coordinate (r2) at (5,2.5);
\coordinate (r3) at (0,4.5);
\coordinate (r4) at (0,.5);
\coordinate (r5) at (5,4);
\coordinate (r6) at (5,1);
\coordinate (r7) at (1.5,5);
\coordinate (r8) at (1.5,0);
\coordinate (r9) at (3,5);
\coordinate (r10) at (3,0);

\fill[color=red!80!white!, opacity=.3] (c) -- (r3) -- (r4) -- (r10) -- (r6) -- cycle;
\fill[color=yellow!80!white!, opacity=.3] (c) -- (r5) -- (r9) -- (r3) -- cycle;
\fill[color=blue!80!white!, opacity=.3] (c) -- (r5) -- (r6) -- (r6) -- (r10) -- (r4) -- cycle;

\path[name path=V, draw] (2.5,2.5) circle [radius=2];
\draw (c) -- (r3);
\draw (c) -- (r4);
\draw (c) -- (r5);
\draw (c) -- (r6);

\node at (2.6,1.2) {\small{$W \cap W^{\omega}$}};
\node at (2.5,3.8) {\small{$L'$}};
\node at (1.1,2.5) {\small{$E$}};
\node at (3.9,2.5) {\small{$F$}};

\end{tikzpicture}

Figure 1

\end{center}
Figure 1 gives a representation of the Witt-Artin decomposition of $V$ with respect to $W$, each of the four slices denoting a direct summand. The circle is all of $V,$ red is used for the subspace $W$ and blue for $W^{\omega},$ giving a violet hue where they intersect. The yellow subspace $L'$ represents a choice of a lagrangian complement of $W\cap W^{\omega}$ in $(E \oplus F)^{\omega}$.

\section{Linear Relations}\label{Linear Relations}

In this section, $W,X,Y$ and $Z$ all denote finite dimensional vector spaces over $\mathbb{R}$. 

\subsection{Definitions, Properties}\label{Definitions, Properties}

 A \textbf{{linear relation}} from $X$ to $Y$ is a linear subspace of the direct sum $X \oplus Y$. In particular, a linear relation is a relation in the set-theoretic sense. If $R \subset X \oplus Y$ is a linear relation, the notation $xRy$ will be used to say that $(x,y) \in R$. 
 
We think of linear relations as generalizations of linear maps in the sense that the graph $\Gamma_F$ of a linear map $F : X \longrightarrow Y$ 
 $$
\Gamma_F = \{(x,y) \in X \oplus Y \mid Fx = y \} \subset X \oplus Y
 $$
is always a linear relation and contains all the information of $F$. 
If $G: Y \longrightarrow Z$ is another linear map, the graph of the composition $G \circ F: X \longrightarrow Z$ is
$$
\Gamma_{GF} = \{(x,z) \in X \oplus Z \mid \exists y\in Y: (x,y) \in \Gamma_F, (y,z) \in \Gamma_G \} \subset X \oplus Z
$$
The composition rule for linear relations should be a generalization of this usual composition rule for maps. Given linear relations $Q \subset X \oplus Y$ and $R \subset Y \oplus Z$ their \textbf{
{composition}} or \textbf{{product}} is defined as 
$$
Q\circ R = \{(x,z) \in X \oplus Z \mid \exists y\in Y: (x,y) \in R, (y,z) \in Q \} \subset X \oplus Z
$$

For a linear relation $R \subset X \oplus Y$ we call $X$ its \textbf{{source}}, $Y$ its \textbf{{target}}, and we define its \textbf{{domain}}, \textbf{{range}},  \textbf{{kernel}} and \textbf{{halo}} respectively as
$$
\text{dom}(R) = \{ x \in X \mid \exists y \in Y : xRy \} 
$$
$$
\text{ran}(R) = \{ y \in Y \mid \exists x \in X : xRy \}
$$
$$
\text{ker}(R) = \{ x \in X \mid xR0 \}
$$
$$
\text{hal}(R) = \{ y \in Y \mid 0Ry \}
$$
A linear relation $R$ is the graph of a linear map if and only if its domain is the entire source space and if it is single valued as a mapping. This is expressed in the conditions
\begin{enumerate}
\item[i)]  $ \forall x \in X \ \exists y \in Y$ such that $xRy$
\item[ii)] $0Rx \ \Rightarrow \ x = 0$
\end{enumerate}
where, by linearity, the second condition is equivalent to saying that if $xRy$ and $xR\tilde{y}$, then  $\tilde{y} = y$. A linear relation is called  \textbf{{cosurjective}} when i) is satisfied, and called \textbf{{coinjective}} if ii) is satisfied. 

The familiar notions of direct sums and the adjoint of a map can be extended to linear relations. If $R \subset X \oplus Y$ and $Q \subset W \oplus Z$ are linear relations, their (external) direct sum is the linear relation in $(X \oplus Y) \oplus (W \oplus Z)$ given by the subspace $R \oplus Q$. The adjoint of a linear map $F: X \longrightarrow Y$ is defined such that for the natural pairing of any vector space $V$ with its dual $V^*$ 
$$
\langle \cdot , \cdot \rangle_V : V^* \times V \rightarrow \mathbb{F} ,\  (l,v) \mapsto l(x)
$$
holds $\langle F^* \alpha,x \rangle_X = \langle \alpha, Fx \rangle_Y $ for all $\alpha \in Y^*, x \in X$,  which, because this pairing is non-degenerate, is equivalent to saying that for any $\alpha \in Y^*, \beta \in X^*$
\begin{equation*}\label{adjoint}
\beta = F^* \alpha \ \Leftrightarrow \ \beta(x) = \alpha(Fx) \quad \forall x \in X
\end{equation*}
One generalizes this to define the adjoint $R^* \subset Y^* \oplus X^*$ of a linear relation $R \subset X \oplus Y$ via the condition  
$$
\alpha R^* \beta \ \Leftrightarrow \ [ \ xRy \Rightarrow \alpha(y) = \beta(x) \ ]
$$
If non-degenerate bilinear maps $B_X : X \times X \rightarrow \mathbb{R}$ and $B_Y : Y \times Y \rightarrow \mathbb{R}$ are given, then 
$$
\tilde B_X : X \rightarrow X^*, \ x \mapsto B_X(x,\cdot) \quad \quad \tilde B_Y : Y \rightarrow Y^*, \ y \mapsto B_Y(y,\cdot)
$$
define natural isomorphisms, which one may use to define the transpose relation $R^t \subset Y \oplus X$ such that the diagram
$$\begin{CD}
X^* @<R^*<<  Y^* \\
@A\tilde B_XAA @AA\tilde B_YA \\
X @<R^t<<  Y
\end{CD}$$
commutes, i.e. 
$$
yR^tx  \  \Leftrightarrow \ [\ zRw \Rightarrow B_Y(y,w) = B_X(x,z) \ ]
$$
When $R$ is a linear map, this is equivalent to 
$$
B_X(R^ty, x) = B_Y(y, Rx) \quad \forall x \in X, y \in Y
$$
i.e. the definition coincides with the usual notion of the transpose of a map. 

One notion which is quite natural in the context of relations is the \textbf{{reverse}} or \textbf{{converse}} $R^r \subset Y \oplus X$ of a linear relation $R \subset X \oplus Y$, which is defined 
$$
(y,x) \in R^r \ \Leftrightarrow \ (x,y) \in R
$$
The operation of taking the converse reverses the roles of source and target; one has
$$ 
\text{dom}(R^r) = \text{ran}(R) \quad \quad \text{ran}(R^r) = \text{dom}(R)
$$  
and $R$ is cosurjective (coinjective) if $R^r$ is surjective (injective). In the special case when $R$ is a linear isomorphism, $R^r$ corresponds to the inverse map $R^{-1}.$ The converse of a relation is {not} in general an inverse though, as illustrated by the simple example when $R = \{(0,y) \mid y \in Y\}$. The composition $R^r R \subset Y \oplus Y$ is equal to all of $Y \oplus Y$, whereas $RR^r \subset X \oplus X$ is equal to $\{(0,0)\}$. The linear relation corresponding to the identity map $X \longrightarrow X$ will be denoted
$$
\Delta_X = \{(x,x) \mid x \in X\}
$$
and this acts as a unit when precomposed with linear relations from $X$ to $Y$ or composed with linear relations from $Y$ to $X$. It may be readily verified that linear relations form the morphisms of a category \textbf{LRel} whose objects are finite dimensional vector spaces, and where for each object $X$, the diagonal $\Delta_X$ is the identity morphism.  

In the following section we will restrict ourselves to considering only linear relations where the source and target space coincide. In this context, we think of the linear relations as the objects, and define a \textbf{morphism} from a linear relation $R \subset X \oplus X$ and a linear relation $Q \subset Y \oplus Y$ as a linear map $S: X \rightarrow Y$ such that 
\begin{equation}\label{morphism}
xRx' \ \Rightarrow (Sx)Q(Sx')
\end{equation}
In this way one again obtains the structure of a category, which we call \textbf{EndLRel}. A morphism $S$ in \textbf{EndLRel} is an \textbf{isomorphism} when it is invertible, and in the special case when $R$ and $Q$ are linear maps, the condition (\ref{morphism}) is equivalent to 
$$
R = S^{-1}QS
$$
i.e. endomorphisms which are isomorphic in \textbf{EndLRel} are precisely those which are conjugate.  

\subsection{Classification}\label{Classification}

We present here the result due to Towber \cite{Towber} stating that any linear relation in  \textbf{EndLRel} is isomorphic to the direct sum of objects in \textbf{EndLRel} which are of only four basic types. Up to the order of summands this decomposition is unique;  it gives
a full classification of the isomorphism classes of \textbf{EndLRel}.  

We describe first the four basic types, using for each dimension $n \in \mathbb{Z}_{\geq 0}$ the model space $V_n = \mathbb{R}e_1 \oplus ... \oplus \mathbb{R}e_n$. The linear endomorphisms of a vector space form a special subset of \textbf{EndLRel} (in fact a subcategory), and by the generalized Jordan normal form any such map is isomorphic to the direct sum of endomorphisms which are indecomposable, i.e. they are not isomorphic to the direct sum of endormorphisms on spaces of smaller dimension. In coordinates, this corresponds to a block matrix form. The decomposition is unique up to order and can be split into two parts, comprising a non-singular and a singular endomorphism respectively. The non-singular part is characterized by the number and size of the blocks, as well as the eigenvalues they correspond to. The singular part can be represented by a sum of Jordan blocks of the form
$$
\left( \begin{array}{ccccc}
0 & 0 & 0 & 0 & 0 \\
1 & 0 & 0 & 0 & 0 \\
0 & 1 & 0 & 0 & 0 \\
\vdots &  & \ddots &  & \vdots \\
0 & 0 & 0 & 1 & 0 
\end{array} \right)
$$
In $V_n \oplus V_n$, for appropriate dimension $n$, such blocks correspond to the linear relation generated over $\mathbb{R}$ by the basis
$$
(e_1,e_2),(e_2,e_3),...,(e_n,0)
$$
Following Towber, we use $\tau^{+} (n)$ to denote this linear relation. 

The other indecomposable basic types identified by Towber have a similar form, but do not correspond to linear endomorphisms (\text{i.e.} they fail either to be single-valued or everywhere-defined). For each dimension $n$ they are denoted $\tau (n)$, $^{+}\tau (n)$ and $^{+}\tau (n)^{+}$, and given in $V_n \oplus V_n$ by the span of 
$$
(e_1,e_2),...,(e_{n-1},e_n),
$$
$$
(0, e_1),(e_1,e_2),..., (e_{n-1},e_n)
$$
and
$$
(0,e_1),(e_2,e_3),...,(e_{n-1},e_n),(e_n,0)
$$
respectively. For completeness we restate\footnote{see \cite{Towber}, p.6} here

\begin{thm}
Every finite-dimensional linear relation is isomorphic to a direct sum of a non-singular linear map and a finite number of linear relations of the types $\tau$, $\tau^{+}$, $^{+}\tau$, and $^{+}\tau^{+}$ (for various values of $n$ and possibly with finite multiplicities). Furthermore the number of summands of each type and each given dimension is unique, i.e. these numbers give a complete set of invariants which classify linear relations up to isomorphism. 
\end{thm}

\section{Linear Canonical Relations}\label{Linear Canonical Relations}

\subsection{Definitions, Properties}\label{Definitions, Properties}

We now come to our main objects of study, in which the structures of linear relations and symplectic linear algebra interact. Let $(X,\omega_X)$ and $(Y,\omega_Y)$ be symplectic vector spaces, and again denote by $Y^-$ the symplectic vector space $(Y,- \omega_Y)$. A \textbf{linear canonical relation} from $X$ to $Y$ is a linear relation $L \subset X \oplus Y^-$ which is a lagrangian subspace with respect to the direct sum symplectic form on $X \oplus Y^-$, i.e. a subspace of dimension $ 1/2(\dim X + \dim Y)$ which is isotropic 
$$
\omega_X(x,x') - \omega_Y(y,y') = 0 \quad \quad \forall \ (x,y), (x',y') \in L
$$
We think of linear canonical relations as a generalization of linear symplectomorphisms (also known as linear canonical transformations). Indeed, if $F: X \rightarrow Y$ is a symplectic map, then $\dim X = \dim Y$ because $F$ is bijective, and by definition
$$
\omega_X(x,x') - \omega_Y(Fx,Fx') = 0
$$
so its graph $\Gamma_{F}$ is an isotropic subspace of $X \oplus Y^-$. Since $\dim \Gamma_{F} = \dim X$ (this holds for any linear map), $\Gamma_{F}$ is lagrangian.

In fact, symplectic maps correspond to the only cases when a linear canonical relation is the graph of a linear map. To see this suppose $L = \Gamma_{f}$ for a linear map $F$ from $X$ to $Y$. Being a graph, $L$ must have dimension equal to $\dim X$, and hence $\dim X = \dim Y$ since $2 \dim L =  \dim X + \dim Y$. So, $F$ is bijective if it is injective. If $(x,0)$ is an element of the kernel of $F$, the condition that $L$ by isotropic gives
$$
\omega_X(x,x') - \omega_Y(0,y') = 0 \quad \forall  y' \in Y, x \in X
$$
which implies that $\omega_X(x,x') = 0\ \forall x' \in X$, so $x=0$ and $F$ is injective. The condition that $\Gamma_{F}$ be isotropic in $X \oplus Y^-$ means that $F$ is symplectic. 

Another important special type of linear canonical relation consists of those which are in some sense the farthest away from being symplectomorphisms. These are linear canonical relations which are entirely either a kernel or a halo, $\text{i.e.}$ of the form
$$
R  = \{ (x,y) \mid  x \in W, y = 0 \} \subset X \oplus Y^- \quad \quad W \subset X \text{ some subspace}
$$
or the converse of a relation of this form. It is easily verified that for such a linear canonical relation one must have $Y = 0$ and $W \subset X$ must be a lagrangian subspace. Thus canonical relations of this form are in one to one correspondence with lagrangian subspaces of $X$ (or subspaces of $Y$ in the case of the converse). 


Similar to linear relations, linear canonical relations are the morphisms of a category, which we call $\textbf{SLRel}$ and where the objects are finite-dimensional symplectic vector spaces. Composition is the same as in $\textbf{LRel}$, and for each symplectic vector space $X$, the identity morphism $1_X$ is again the diagonal $\Delta_X = \{(x,x) \mid x \in X\}$ (which, as the graph of the identity, is indeed a linear canonical relation). To show that the composition
$$
L_2 \circ L_1 = \{ (x,z) \mid \exists y \in Y : (x,y) \in L_1, (y,z) \in L_2 \}
$$
of linear canonical relations $L_1 \subset X \oplus Y^-$, $L_2 \subset Y \oplus Z^-$ is again a linear canonical relation, we describe the subspace $L_2 \circ L_1$ as the image of a lagrangian under a reduction map. To see this, we note first that the set $L_2 \circ L_1$ is the result of the following steps: 
\begin{enumerate}
\item[i)] intersect $L_1 \oplus L_2 \subset X \oplus Y^- \oplus Y \oplus Z^-$ with $X \oplus \Delta_Y \oplus Z^-$
\item[ii)] project $(L_1 \oplus L_2) \cap (X \oplus \Delta_Y \oplus Z)$ onto $X \oplus Z$ 
\end{enumerate}
It is easily checked that $C = X \oplus \Delta_Y \oplus Z$ is a coisotropic subspace of $X \oplus Y^- \oplus Y \oplus Z^-$ (its orthogonal space is $0 \oplus \Delta_Y \oplus 0$) and that $L_1 \oplus L_2$ is a lagrangian subspace in $X \oplus Y^- \oplus Y \oplus Z^-$. Hence, by Lemma \ref{coisotropic reduction} the reduction map $C \rightarrow C/(0 \oplus \Delta_Y \oplus 0)$ maps $(L_1 \oplus L_2) \cap (X \oplus \Delta_Y \oplus Z)$ to a lagrangian subspace. But $C/(0 \oplus \Delta_Y \oplus 0) = X \oplus 0 \oplus 0 \oplus Z^-$, which one identifies with $X \oplus Z^-$, and the image of $(L_1 \oplus L_2) \cap (X \oplus \Delta_Y \oplus Z)$ under the reduction and after this identification is precisely the image of the projection in step ii) above, i.e. $L_2 \circ L_1$.

In general, the reduction of a symplectic space $X$ by a coisotropic subspace $C$, given by the reduction map $\rho: C \longrightarrow C/C^{\omega}$, can be recast in the present context as a linear relation $R \subset X \oplus {C/C^{\omega}}^-$ which is surjective, single-valued (coinjective), but whose domain is $C$, \text{i.e.} it is not defined everywhere. The fact that $\rho^*[\omega_X] = \omega_X$ means precisely that $R$ in this case is an isotropic subspace of $X \oplus C/{C^{\omega}}^-$ and 
$$2 \dim R = 2 \dim C = \dim(C/C^{\omega}) + \dim C^{\omega} + \dim C = \dim(C/C^{\omega}) + \dim V$$ 
shows that $R$ is lagrangian. It turns out that any linear canonical relation which is surjective and coinjective is induced in the above way by a reduction map on some coisotropic subspace. For this reason one calls a surjective, coinjective canonical relation a \textbf{{reduction}} and accordingly a cosurjective, injective canonical relation a \textbf{{coreduction}}. Equivalently, a coreduction is simply the converse of a reduction. 

It is worth noting that in the literature one usually refers here to the transpose instead of the converse. Either wording is appropriate since for linear canonical relations these two concepts coincide with respect to $B_X = \omega_X$ on each symplectic vector space $X$.  Indeed, if $R \subset X \oplus Y^-$ is a canonical relation and $\bar \omega = \omega_{X \oplus Y^-}$, then
\begin{align*}
 yR^tx  \  & \Leftrightarrow \ [\ zRw \Rightarrow B_X(x,z) = B_Y(y,w) \ ] \\
		& \Leftrightarrow [\ zRw \Rightarrow \omega_X(x,z) = \omega_Y(y,w) \ ] \\
		& \Leftrightarrow [\ zRw \Rightarrow \bar \omega((y,x),(z,w))=0  \ ] \\
		& \Leftrightarrow \ (x,y) \in R^{\bar \omega} = R \\
		& \Leftrightarrow \ xRy
\end{align*}

\subsection{Factorization}\label{Factorization}

In this section we show how every linear canonical relation may be factored as the composition of a reduction, a symplectic map and a coreduction. This has in particular important consequences for our classification problem. We once again restrict ourselves to the special context of linear canonical relations whose source and target coincide, fixing our notation such that  $L \subset V \oplus V^-$ always denotes a linear canonical relation, $A$ denotes its domain, $B$ its image, and $\bar \omega$ the symplectic form on $V \oplus V^-$. 

If $\hat L \subset \hat V \oplus \hat V^-$ is another canonical relation and $S: V \rightarrow \hat V$ an equivalence between the two, \text{i.e.} a symplectic map such that 
$$
xLy \Leftrightarrow (Sx)\hat L(Sy)
$$
then it is clear that $S$ maps $A$ and $B$ respectively to the domain and range of $\hat L$,
\begin{equation}\label{coisotropic pair condition}
S(A) = \hat A, \quad  S(B) = \hat B
\end{equation}
Another simple but important observation is that the orthogonal of $A$ in $V$ is the kernel of $L$. An analogous statement holds for $B$ and the halo of $L$. 

\begin{lem} \label{A^perp is zero related} For the domain $A$ of  $L \subset V \oplus V^-$holds: $v \in V$ is in $A^{\omega}$ iff $(v,0) \in L$. 
\end{lem}

\pf Assume first that $v \in A^{\omega}$. For every $(a,b) \in L$ one has $\overline{\omega}((v,0),(a,b)) =  \omega(v,a) - \omega(0,b) = 0$, and hence $(v,0) \in L^{\omega}$, and $L^\omega = L$ because $L$ is lagrangian. Assume on the other hand that $(v,0) \in L$. Then $\omega(v,a) = \overline{\omega}((v,0),(a,b)) - \omega(0,b) = 0$ for every $(a,b) \in L$, and hence $\omega(v,a) = 0$ for every $a \in A$, i.e. $v\in A^{\omega}$. 
\qed
This in turn leads to the following factorization result\footnote{see \cite{Benenti Tulczyjew}, Proposizioni 4.4 \& 4.5 }.

\begin{prop}[Benenti-Tulczyjew] \label{Tulczyjew} For any $L \subset V \oplus V^-$, the domain $A$ and range $B$ are coisotropic subspaces of $V$ and the quotient relation $[L] \subset A/A^{\omega} \oplus B/B^{\omega}$ induced by $L$ defines a symplectomorphism between the reduced spaces $A/A^{\omega}$ and $B/B^{\omega}$. In other words $L$ factors as $R_B^t \circ [L] \circ R_A$, where $R_A$ is the canonical relation of the reduction of $V$ by $A$, and $R_B^t$ the transpose of the reduction of $V$ by $B$. 
\end{prop}

\pf $A$ (and $B$) must be coisotropic, i.e. $A^{\omega} \subset A$, since by the previous proposition one has  $(v,0) \in L$ for any $v\in A^{\omega}$, and hence also $v\in A$.  We recall that the induced quotient relation $[L]\in A/A^{\omega} \times B/B^{\omega}$ is (well)defined by $([a],[b]) \in [L]$ iff $(a,b) \in L$. Observe that for every $[a] \in A/A^{\omega} $ there exists $[b] \in B/B^{\omega}$ such that $([a],[b]) \in L$, simply by virtue of the fact that $a \in A $ implies that $(a,b) \in L$ for some $b \in B$. Thus $[L]$ is the graph of a linear map $A/A^{\omega} \rightarrow B/B^{\omega}$. Via the same argument, but with the roles of $a$ and $b$ reversed, one sees that this linear map must also be surjective. For injectivity, assume that $([a],[b])$ and $([\tilde{a}],[b])$ are both in $[L]$. By linearity, $([a]-[\tilde{a}],[0])$ is also in $[L]$, and hence $(\tilde{a} - a,0) \in L $, which by proposition~\ref{A^perp is zero related}  is equivalent to $\tilde{a} - a \in A^{\omega}$. It follows that $[\tilde{a} - a] = [0]$, i.e. $[\tilde{a}] = [a]$. So, the linear mapping corresponding to $[L]$ is nonsingular. To check that it is a symplectomorphism, let $([a],[b])$ and $([c],[d])$ be two elements of $[L]$. If $[\omega]_A$ and $[\omega]_B$ denote the induced symplectic forms on $A/A^{\omega}$ and $B/B^{\omega}$ respectively, one has $[\omega]_B ([b],[d]) - [\omega]_A ([a],[c]) = \omega(b,d) - \omega(a,c) = \overline{\omega}((a,b),(c,d))$, which is zero as desired because $(a,b)$ and $(c,d)$ are in $L$. 
\qed

\begin{cor} \label{symplectic parts are related}
Let $A_1$ and $B_1$ be such that $A = A^{\omega} \oplus A_1$ and $B = B^{\omega} \oplus B_1$. Denote by $u = u_0 + u_1$ the corresponding decomposition of a vector $u$ in $A$ or $B$. Then $(v,w) \in L \Rightarrow (v_1, w_1) \in L$.
\end{cor}

\pf
One has $(v,w) =(v_0 + v_1, w_0 + w_1) = (v_1, w_1) +  (v_0,0) + (0, w_0)$, or equivalently $(v_1, w_1) = (v,w) - (v_0,0) - (0, w_0)$. The terms $(v_0,0)$ and $(0, w_0)$ are in $L$ by Proposition \ref{A^perp is zero related}, hence if $(v,w)$ is in $L$, then so is $(v_1, w_1)$. 
\qed

\begin{cor} One has $dim(A^{\omega}) = dim(B^{\omega})$, and hence also $dimA = dim B$. 
\end{cor}

\pf Let $A_1$ be a subspace of $A$ such that $A = A^{\omega} \oplus A_1$, and define a $B_1$ analogously. Because $A_1 \simeq A/A ^{\omega}$ and $B_1 \simeq B/B ^{\omega}$, by Proposition~\ref{Tulczyjew} it follows that $\dim A_1 = \dim B_1$. Also, we know that $2n = \dim A + \dim A^{\omega} = \dim A_1 + 2 \dim (A^{\omega})$, and similarly so for $B$. Combining these facts gives the result. 
\qed



\begin{cor}
Given $L \subset V \oplus V^-$ with domain $A$ and range $B$, let $A_1$ and $B_1$ be any choice of subspaces such that $A = A^{\omega} \oplus A_1$, $B = A^{\omega} \oplus B_1$. Then $L$ induces a symplectic map $\phi_L : A_1 \rightarrow B_1$ such that  $\forall \ (v_1,w_1) \in A_1 \times B_1: (v_1,w_1) \in L \ \Leftrightarrow \ w_1 = \phi_L (v_1)$. 
\end{cor}

\begin{rem}
The notation $\phi_L$ does not reflect the fact that this map depends on the choice of $A_1$ and $B_1$. 
\end{rem}

\begin{prop} \label{equivalence by parts}
Let $L$ and $\hat L$ be two canonical relations, $A,B$ and $\hat A, \hat B$ their respective domains and ranges, and $A_1$ and $B_1$ any subspaces such that $A = A^{\omega} \oplus A_1$, $B = B^{\omega} \oplus B_1$. Let $\phi_L:A_1 \rightarrow B_1$ be the symplectic map induced by $L$ and these decompositions. A symplectic map $S:V \rightarrow V$ is an equivalence between $L$ and $\hat L$ if and only if
\begin{enumerate}
\item[i)] $S(A^{\omega}) = \hat A^{\omega}$ and $S(B^{\omega}) = \hat B^{\omega}$
\item[ii)] $ \phi_{\hat L} \circ S \vert_{A_1} = S \vert_{B_1} \circ \phi_L $
\end{enumerate}
whereby $\phi_{\hat L}: S(A_1) \rightarrow S(B_1)$ is the symplectic map induced by $\hat L$ and the decompositions $\hat A = \hat A^{\omega} \oplus S(A_1)$ and $\hat B = \hat B^{\omega} \oplus S(B_1)$. 
\end{prop}

\pf
Assume i) and ii) hold. We show first $(v,w) \in L \Rightarrow (Sv,Sw) \in \hat L$. Let $(v,w) \in L$. We have $(Sv,Sw) =(Sv_0 + Sv_1, Sw_0 + Sw_1) = (Sv_1, Sw_1) +  (Sv_0,0) + (0, Sw_0)$, where $v_0 \in A^{\omega}$, $w_0 \in B^{\omega}$ and $(v_1, w_1) \in A_1 \times B_1$. Since $Sv_0 \in \hat A ^{\omega}$, we have $(Sv_0,0) \in \hat L$, analogously $(0, Sw_0) \in \hat L$, and $(Sv_1, Sw_1) \in \hat L$ follows from assumption ii): $(Sv_1,Sw_1) \in \hat L \Leftrightarrow S\vert_{B_1} w_1 = \phi_{\hat L }(S\vert_{A_1} v_1) \Leftrightarrow S\vert_{B_1} w_1 = S\vert_{B_1} (\phi_L (v_1)) \Leftrightarrow w_1 = \phi_L (v_1) \Leftrightarrow (v_1,w_1) \in L$, and indeed $(v_1,w_1) \in L$ follows from $(v,w) \in L$ (cf. Corollary \ref{symplectic parts are related}). So all three summands are in $\hat L$, and hence so is their sum $(Sv,Sw)$. Because $S$ is invertible and symplectic, the converse implication follows by arguing symmetrically in the opposite direction. 

Now assume $S$ is an equivalence between $L$ and $\hat L$. The property i) follows via Proposition \ref{A^perp is zero related}, and implies that $(v_1,w_1) \in A_1 \times B_1$ iff $(v,w) \in L$ and $(Sv_1,Sw_1) \in S(A_1) \times S(B_1)$ iff $(Sv,Sw) \in L$, and hence we have $(v_1,w_1) \in L$ iff $(Sv_1,Sw_1) \in \hat L$ . Property (2) then follows from the set of equivalences $ S\vert_{B_1} w_1 = \phi_{\hat L }(S\vert_{A_1} v_1) \Leftrightarrow (Sv_1,Sw_1) \in L \Leftrightarrow (v_1,w_1) \in L \Leftrightarrow w_1 = \phi_L (v_1)$, which hold for all $(v_1,w_1)  \in A_1 \times B_1$.
\qed
This last proposition breaks our classification problem into two parts. The property i) above is equivalent to (\ref{coisotropic pair condition}), \text{i.e.} the condition
$$
S(A) = \hat A, \quad S(B)= \hat B
$$
and constitutes a necessary step in the classification of linear canonical relations. For this reason we now investigate the question of when such a symplectic map $S$ exists between any two given pairs $(A,B)$ and $(\hat A, \hat B)$ of coisotropic subspaces. 

\section{Coisotropic Pairs}\label{Coisotropic Pairs}

The results in this section constitute joint work together with Alan Weinstein which have, in the meantime, been extended to the settings of presymplectic and Poisson vector spaces (see \cite{Lorand Weinstein}).

We call an ordered pair $(A,B)$ of coisotropic subspaces $A,B \subset V$ a \textbf{{coisotropic pair}} and say that coisotropic pairs $(A,B)$ and $(\hat A, \hat B)$ given in $(V,\omega)$ and $(\hat V, \hat \omega)$ respectively are \textbf{{equivalent}} if there exists a linear symplectomorphism $S:V \rightarrow \hat V$ such that $S(A) = \hat A$ and $S(B) = \hat B$. For a coisotropic pair $(A,B)$ in $(V,\omega)$ we allow the general situation where $\dim A$ and $\dim B$ are not necessarily equal; we will see that $(A,B)$ is fully characterized up to equivalence by the following five numbers
$$\dim (A^{\omega} \cap B^{\omega}), \ \ \dim A^{\omega}, \ \  \dim B^{\omega}, \ \  \frac{1}{2}\dim V, \ \ \dim (A^{\omega} \cap B)$$  
which will be called the \textbf{{canonical invariants}} of $(A,B)$ and labeled $k_1$ through $k_5$ in the above order. They are largely independent, subject only to certain inequalities (see Corollary \ref{ineq only constraint}). 
 
The first four invariants $k_1, k_2, k_3, k_4$ characterize the subspaces $A$ and $B$ up to the above equivalence if one drops the condition that $S$ be symplectic and that $A$ and $B$ be coisotropic, i.e. these four invariants contain the purely linear algebraic information.  Indeed, using the identities $\dim W^{\omega} = \dim V - \dim W$ and $(E + F)^{\omega} = E^{\omega} \cap F^{\omega}$, which hold for any subspaces $W$, $E$, $F$ in $V$, one can obtain the the linear algebraic data 
\begin{equation}\label{lin alg invariants}
\dim V, \ \dim A, \ \dim B, \ \dim (A \cap B)
\end{equation}
from these four invariants:
$\dim V = 2 \cdot \frac{1}{2} \dim V$, $\dim A = \dim V - \dim A^{\omega}$, $\dim B = \dim V - \dim B^{\omega}$,
and
\begin{align*}
\dim (A \cap B) & = \dim V - \dim (A \cap B)^{\omega} \\
& = \dim V - \dim (A^{\omega} + B^{\omega}) \\
 & = \dim V - \dim A^{\omega} - \dim B^{\omega} + \dim (A^{\omega} \cap B^{\omega})
\end{align*}
It is straightforward to check that his relationship is invertible; one could thus equivalently use the numbers (\ref{lin alg invariants}) as the first four invariants. 

The fifth invariant $k_5 = \dim (A^{\omega} \cap B)$ is what fixes the symplectic information. One could equivalently choose $ \dim (B^{\omega} \cap A)$ as the fifth invariant, since
\begin{align*}
 \dim (A^{\omega} \cap B) & = \dim V - \dim (A^{\omega} \cap B)^{\omega} \\
  & = \dim V - [\dim A + \dim B^{\omega} -  \dim (B^{\omega} \cap A)]  \\
  & = \dim (B^{\omega} \cap A) + \dim A^{\omega} -\dim B^{\omega}
\end{align*} 

When $\dim A = \dim B$ it follows that $\dim (A^{\omega} \cap B) = \dim (B^{\omega} \cap A)$, and a total of four invariants suffice to characterize the coisotropics $A$ and $B$. They can be given in a symmetric way as 
$$
\dim V, \ \dim (A + B), \ \dim (A \cap B), \ \text{rank} (A^{\omega} + B^{\omega})
$$
where for any subspace $W \subset V$, rank$(W) = \dim W - \dim (W \cap W^{\omega})$. The symmetry of these invariants implies that $(A,B)$ and $(B,A)$ are equivalent as coisotropic pairs in this special case. 

Note that because a coisotropic subspace $A$ is uniquely determined by the isotropic subspace $A^{\omega}$, and $S(A^{\omega}) = S(A)^{\hat \omega}$ for any linear symplectomorphism $S:V \rightarrow \hat V$, one could equivalently consider isotropic pairs instead of coisotropic ones. This indeed simplifies some calculations and proofs; for the present though we treat things from the coisotropic standpoint.

\subsection{General classification of coisotropic pairs}\label{General classification of coisotropic pairs}

Given a coisotropic pair $(A,B)$, we fix the notation $I := A^{\omega} \cap B^{\omega} $ and $K:= A^{\omega} \cap B + B^{\omega} \cap A$. As announced, the numbers $\dim (A^{\omega} \cap B^{\omega})$, $\dim A^{\omega}$, $\dim B^{\omega}$, $ 1/2 \dim V$ and $\dim (A^{\omega} \cap B)$, which we call the canonical invariants associated to $(A,B)$, completely characterize a coisotropic pair up to equivalence.

\begin{prop}\label{main prop}
Let $(A,B)$ and $(\hat A,\hat B)$ be pairs of coisotropic subspaces in $(V,\omega)$ and $(\hat V, \hat \omega)$ respectively. Then $(A,B)$ and $(\hat A,\hat B)$ are equivalent if and only if their associated canonical invariants are equal. 
\end{prop}

\pf
If $(A,B)$ and $(\hat A,\hat B)$ are equivalent via some symplectic map $S: V \rightarrow \hat V$, it is clear that all the canonical invariants of $(A,B)$ and $(\hat A,\hat B)$ coincide. 

For the converse, we will show that $V$ can be written as an $\omega$-orthogonal direct sum of five symplectic subspaces
$$ V = D \overset{\omega}{\oplus} E \overset{\omega}{\oplus} F \overset{\omega}{\oplus} G \overset{\omega}{\oplus} H $$
where each symplectic piece, except for $F$, is further decomposed as a lagrangian pair
$$  D = I \oplus J, \ \ \  E = E_1 \oplus E_2 , \ \ \  G = G_1 \oplus G_2, \ \ \  H = H_1 \oplus H_2 $$
so that we obtain a decomposition of $V$ into a total of nine subspaces
\begin{equation}\label{inner decomp}
V = (I \oplus J) \overset{\omega}{\oplus} (E_1 \oplus E_2) \overset{\omega}{\oplus} F \overset{\omega}{\oplus} (G_1 \oplus G_2) \overset{\omega}{\oplus} (H_1 \oplus H_2)
\end{equation}
Moreover, this decomposition will have the following properties:

\begin{enumerate}[i)]
\item  the dimension of each summand is uniquely determined by the canonical invariants of $(A,B)$
\item $A$ and $B$ are decomposable as
$$A = I \oplus E_1 \oplus G_1 \oplus F \oplus H_1 \oplus H_2$$ 
$$ B = I \oplus E_2  \oplus H_1 \oplus F  \oplus G_1 \oplus G_2$$ 

\end{enumerate} 

One can decompose $\hat V$ in an analogous manner, and hence when $(A,B)$ and $(\hat A ,\hat B)$ have the same invariants, by property i) the dimensions of corresponding symplectic pieces in the decompositions of $V$ and $\hat V$ will match. In this case, for dimension reasons alone there exist five symplectic maps, one each between corresponding symplectic pieces, i.e. one from $D$ to $\hat D$, one from $E$ to $\hat E$, and so on.  These maps can further be be chosen to respect the respective decompositions into lagrangian pairs. 

Because the five-part decompositions of $V$ and $\hat V$ are symplectic-orthogonal, the direct sum of these five symplectic maps defines a symplectic map $S: V \rightarrow \hat V$ which respects all nine summands of the decompositions of $V$ and $\hat V$. In particular, by property ii), $S$ will then also satisfy $S(A)=\hat A$ and $S(B)=\hat B$.

To achieve the decomposition (\ref{inner decomp}) we will construct a certain Witt-Artin decomposition of $V$ with respect to $W := A^{\omega} + B^{\omega}$, refined and adapted to the coisotropic subspaces $A$ and $B$.

Recall that $I = A^{\omega} \cap B^{\omega} $ and $K= A^{\omega} \cap B + B^{\omega} \cap A$, and note that $$W^{\omega} = (A^{\omega} + B^{\omega})^{\omega} = A \cap B$$ 
and
$$W \cap  W^{\omega} = (A^{\omega} + B^{\omega}) \cap (A \cap B) = A^{\omega} \cap B + B^{\omega} \cap A = K$$ 
We begin by decomposing $A^{\omega}$ into three parts by choosing a subspace $G_1$ such that $A^{\omega} \cap B = I \oplus G_1$ and a subspace $E_1$ such that $A^{\omega} = A^{\omega} \cap B \oplus E_1$, giving a decomposition 
$$ A^{\omega} = I \oplus G_1 \oplus E_1 $$
Analogously we obtain a decomposition 
$$B^{\omega} = I \oplus H_1 \oplus E_2$$
where $E_2$ is such that $B^{\omega} = B^{\omega} \cap A \oplus E_2$, and $H_1$ such that $B^{\omega} \cap A = I \oplus H_1$. Note that $H_1$ and $G_1$ have zero intersection, since $H_1 \cap G_1 \subset A^{\omega} \cap B^{\omega} = I$ and $H_1 \cap I = 0$ and $G_1\cap I = 0$. Similarly, $E_1 \cap E_2 = 0$. In particular we have 
$$K = A^{\omega} \cap B + B^{\omega} \cap A =  I \oplus G_1 + I \oplus H_1 =  I \oplus G_1 \oplus H_1$$ 

We now set $E:= E_1 \oplus E_2$. This defines a subspace such that $K \oplus E = A^{\omega} + B^{\omega} = W$. Indeed, 
$$K + E = I \oplus G_1 \oplus H_1 + E_1 \oplus E_2 = A^{\omega} + B^{\omega} = W$$ and $K \cap E = 0$ since 
\begin{align*}
\dim K + \dim E & = (\dim I  + \dim G_1 + \dim H_1) + (\dim E_1 + \dim E_2) \\
& = \dim A^{\omega} + \dim B^{\omega} - \dim I \\
& = \dim (A^{\omega} + B^{\omega}) = \dim W \\
& = \dim (K + E) 
\end{align*}
Because $E$ is a complement of $K =W\cap W^{\omega}$ in $W$, $E$ is symplectic by Lemma \ref{Witt-Artin decomposition}, and since $E_1$ and $E_2$ are both isotropic, we conclude that they form a transversal lagrangian pair in $E$. 

To obtain a Witt-Artin decomposition with respect to $W$, we choose a complement $F$ of $W \cap W^{\omega} = K$ in $W^{\omega} = A \cap B$, i.e. so that 
$$ A \cap B  =  K \oplus F $$
Applying Lemma \ref{Witt-Artin decomposition} again we know that $F$ is symplectic, as is $E$, and $V$ decomposes into the $\omega$-orthogonal direct sum
$$ V = E \overset{\omega}{\oplus} F  \overset{\omega}{\oplus} (E \oplus F)^{\omega} $$
with $K$ as a lagrangian subspace of the symplectic subspace $(E \oplus F)^{\omega}$.

We refine this decomposition by choosing a lagrangian complement $K'$ of $K$ in $(E \oplus F)^{\omega}$ and by defining a decomposition in $K'$ using the decomposition $K =  I \oplus G_1 \oplus H_1$ as follows. Any basis $\textbf{q}$ of $K$ is mapped under $\tilde{\omega}$ to a basis of $(K')^*$, whose dual basis $\textbf{p}$ in $K'$ is conjugate to $\textbf{q}$, \text{i.e.} together $\textbf{q}$ and $\textbf{p}$ form a symplectic basis of $K \oplus K'$.  If we consider a basis $\textbf{q}$ which is adapted to the decomposition in $K$, then this partitioning induces a partitioning of $\textbf{p}$ which defines subspaces $J$, $G_2$ and $H_2$ in $K'$ such that 
$$
K' = J \oplus G_2 \oplus H_2
$$
and $D:= I \oplus J$, $G:= G_1 \oplus G_2$ and $H:=H_1\oplus H_2$ are $\omega$-orthogonal symplectic subspaces, comprised each of a lagrangian pair, giving
$$
K \oplus K' = D \oplus G \oplus H
$$
In total we thus obtain a decomposition 
$$ V = (I \oplus J ) \oplus (E_1 \oplus E_1) \oplus F \oplus (G_1 \oplus G_2) \oplus (H_1 \oplus H_2) $$ 
where parentheses enclose transversal lagrangian pairs in a symplectic subspace. This decomposition is visualized in Figure 2 - the full circle represents $V$, each piece is a direct summand, and lagrangian pairs are aligned symmetrically with respect to the horizontal axis and shaded with colors of a similar hue. 

\begin{center}
\begin{tikzpicture}[scale = 1]

\path[name path=V, draw] (2.5,2.5) circle [radius=2];

\path[name path=V, clip] (2.5,2.5) circle [radius=2];

\coordinate (c) at (2.5,2.5);
\coordinate (r1) at (0,2.5);
\coordinate (r2) at (5,2.5);
\coordinate (r3) at (0,4.5);
\coordinate (r4) at (0,.5);
\coordinate (r5) at (5,4);
\coordinate (r6) at (5,1);
\coordinate (r7) at (1.5,5);
\coordinate (r8) at (1.5,0);
\coordinate (r9) at (3,5);
\coordinate (r10) at (3,0);

\fill[color=brown!80!white!, opacity=.4] (c) -- (r3) -- (r1) -- cycle; 
\fill[color=brown!80!gray!, opacity=.6] (c) -- (r4) -- (r1) -- cycle; 
\fill[color=yellow!80!white!, opacity=.4] (c) -- (r7) -- (r3) -- cycle; 
\fill[color=yellow!80!gray!, opacity=.6] (c) -- (r8) -- (r4) -- cycle; 
\fill[color=orange!80!white!, opacity=.4] (c) -- (r9) -- (r7) -- cycle; 
\fill[color=orange!80!gray!, opacity=.6] (c) -- (r10) -- (r8) -- cycle; 
\fill[color=red!80!white!, opacity=.4] (c) -- (r5) -- (r9) -- cycle; 
\fill[color=red!80!gray!, opacity=.6] (c) -- (r6) -- (r10) -- cycle; 
\fill[color=violet!80!brown!, opacity=.5] (c) -- (r5) -- (r6) -- cycle; 

\path[name path=V, draw] (2.5,2.5) circle [radius=2];
\draw (c) -- (r1);
\draw (c) -- (r3);
\draw (c) -- (r4);
\draw (c) -- (r5);
\draw (c) -- (r6);
\draw (c) -- (r7);
\draw (c) -- (r8);
\draw (c) -- (r9);
\draw (c) -- (r10);

\node at (2.35,1.2) {\small{$I$}};
\node at (2.35,3.8) {\small{$J$}};
\node at (1.2,2) {\small{$E_1$}};
\node at (1.2,2.9) {\small{$E_2$}};
\node at (1.7,1.4) {\small{$G_1$}};
\node at (3.3,1.3) {\small{$H_1$}};
\node at (3.9,2.5) {\small{$F$}};
\node at (1.7,3.6) {\small{$G_2$}};
\node at (3.3,3.6) {\small{$H_2$}};

\end{tikzpicture}

Figure 2

\end{center}
\ 

The coisotropics $A$ and $B$ are related to the decomposition in $K'$ in that $G_2 = B \cap K'$ and $H_2 = A \cap K'$. To see this it suffices to show that their $\omega$-orthogonal complements are equal. For the case of $A\cap K'$ (the case for $B \cap K'$ is analogous) one has
\begin{align*}
(A \cap K' )^{\omega} & = A^{\omega} + (K')^{\omega} \\
				& = I \oplus G_1 \oplus E_1 + E \oplus F \oplus K' \\
				& = I \oplus G_1 \oplus E \oplus F \oplus K' \\
				& = H_2^{\omega}
\end{align*}
where we use in the last step that $H_2$ is $\omega$-orthogonal to $D$, $G$, $K'$ and $E \oplus F$ and that the dimensions match. 
\

It can now be quickly checked that our decomposition of $V$ satisfies property ii), \text{i.e.} that
$$A = I \oplus E_1 \oplus G_1 \oplus F \oplus H_1 \oplus H_2$$ 
$$ B = I \oplus E_2  \oplus H_1 \oplus F  \oplus G_1 \oplus G_2$$ 
We show this for $A$, the decomposition of $B$ follows in the same way. The inclusion ``$\supset$" is obvious since all the spaces on the right-hand side are subsets of $A$. The opposite inclusion ``$\subset$" can be argued using dimensions:
\begin{align*}
\dim A & =  \dim V - \dim A^{\omega} \\
	       & =  \dim(I \oplus E_1 \oplus G_1  \oplus F \oplus H_1 \oplus H_2) + \dim (J \oplus G_2 \oplus E_2) - \dim A^{\omega} \\
	        & =  \dim(I  \oplus E_1 \oplus G_1  \oplus F \oplus H_1 \oplus H_2) 
\end{align*}
where the last equality follows from the fact that 
$$ \dim A^{\omega} = \dim (I  \oplus E_1 \oplus G_1 ) = \dim (J \oplus E_2 \oplus G_2) $$
since $\dim I = \dim J$, $ \dim E_1 = \dim E_2$, and $\dim G_1 = \dim G_2$ (each pair of subspaces is a lagrangian pair in $D$, $E$ and $G$ respectively). 

The decompositions of $A$ and $B$ are visualized below.  
\begin{center}
\begin{minipage}{.45\textwidth}
\begin{center}

\begin{tikzpicture}[scale = .95]

\path[name path=V, draw] (2.5,2.5) circle [radius=2];

\path[name path=V, clip] (2.5,2.5) circle [radius=2];

\coordinate (c) at (2.5,2.5);
\coordinate (r1) at (0,2.5);
\coordinate (r2) at (5,2.5);
\coordinate (r3) at (0,4.5);
\coordinate (r4) at (0,.5);
\coordinate (r5) at (5,4);
\coordinate (r6) at (5,1);
\coordinate (r7) at (1.5,5);
\coordinate (r8) at (1.5,0);
\coordinate (r9) at (3,5);
\coordinate (r10) at (3,0);

\fill[color = yellow!70!white!, opacity = .3] (2.5,2.5) circle [radius=2]; 

\fill[color=blue!50!white!, opacity=.5] (c) -- (r4) -- (0,0) -- (5,0) -- (r5) -- cycle; 
\fill[color=blue!50!white!, opacity=.5] (c) -- (r7) -- (r3) -- cycle; 
\fill[color=blue!50!white!, opacity=.5] (c) -- (r3) -- (r1) -- cycle; 
\fill[color=green!50!white!, opacity=.5] (c) -- (r1) -- (0,0) -- (5,0) -- (5,5) -- (3,5) -- cycle; 

%
%
%

\fill[color=cyan!70!white!, opacity=.3] (c) -- (r3) -- (r1) -- cycle; 
\fill[color=green!70!yellow!, opacity=.2] (c) -- (r4) -- (r1) -- cycle; 

\fill[color=green!70!yellow!, opacity=.2] (c) -- (r8) -- (r4) -- cycle; 
\fill[color=green!70!yellow!, opacity=.2] (c) -- (r10) -- (r8) -- cycle; 
\fill[color=cyan!70!white!, opacity=.3] (c) -- (r10) -- (r8) -- cycle; 
\fill[color=cyan!70!white!, opacity=.3] (c) -- (r6) -- (r10) -- cycle; 




\path[name path=V, draw] (2.5,2.5) circle [radius=2];
\draw (c) -- (r1);
\draw (c) -- (r3);
\draw (c) -- (r4);
\draw (c) -- (r5);
\draw (c) -- (r6);
\draw (c) -- (r7);
\draw (c) -- (r8);
\draw (c) -- (r9);
\draw (c) -- (r10);

\node at (2.35,1.2) {\small{$I$}};
\node at (2.35,3.8) {\small{$J$}};
\node at (1.2,2) {\small{$E_1$}};
\node at (1.2,2.9) {\small{$E_2$}};
\node at (1.7,1.4) {\small{$G_1$}};
\node at (3.3,1.3) {\small{$H_1$}};
\node at (3.9,2.5) {\small{$F$}};
\node at (1.7,3.6) {\small{$G_2$}};
\node at (3.3,3.6) {\small{$H_2$}};

\end{tikzpicture}

Figure 3

\end{center}
\end{minipage}%
\begin{minipage}{.45\textwidth}
\begin{center}

\begin{tikzpicture}[scale = .7]

\path[use as bounding box] (0,0) rectangle (8, 5.2);

\draw[name path=border] (0,0) rectangle (8,5.2);
\fill [name path=border, color=yellow!70!white!, opacity=.3] (0,0) rectangle (8,5.2);

\path[name path global=B] (3.5,2.5) ellipse [x radius=3, y radius=1, rotate=40]; 
\path[name path global=A] (4,2.5) ellipse [x radius=4, y radius=1.5, rotate=155];
\path[name path global=B'] (4,2) ellipse [x radius=4, y radius=1, rotate=40]; 
\path[name path global=A'] (3.5,2) ellipse [x radius=5, y radius=1.5, rotate=155]; 


\path[name path global=B, fill] [color=blue!50!white!, opacity=.5] (3.5,2.5) ellipse [x radius=3, y radius=1, rotate=40]; 
\path[name path global=A, fill] [color=green!50!white!, opacity=.5] (4,2.5) ellipse [x radius=4, y radius=1.5, rotate=155];
\begin{scope}
\path[name path global=B, clip] (3.5,2.5) ellipse [x radius=3, y radius=1, rotate=40];
\fill [name intersections={of=A' and B, name=q, total=\t}, color=cyan!70!white!, opacity=.3]
(q-4) -- (3,5.2) -- (8,5.2) -- (8,3) -- (q-3) -- cycle;
\end{scope}

\begin{scope}
\path[name path global=A, clip] (4,2.5) ellipse [x radius=4, y radius=1.5, rotate=155]; 
\fill [name intersections={of=A and B', name=p, total=\t}, color=green!70!yellow!, opacity=.2]
(p-1) -- (0,2) -- (0,5.2) -- (5,5.2) -- (p-3) -- cycle;
\end{scope}

\path[name path global=B, draw] (3.5,2.5) ellipse [x radius=3, y radius=1, rotate=40]; 
\path[name path global=A, draw] (4,2.5) ellipse [x radius=4, y radius=1.5, rotate=155];

\draw [name intersections={of=A' and B, name=q, total=\t}]
(q-4) -- (q-3);
\draw [name intersections={of=A and B', name=p, total=\t}]
(p-1) -- (p-3);

\node at (4.1,3.7) {\small{$I$}};
\node at (7,3.5) {\small{$J$}};
\node at (1.7,3.8) {\small{$E_1$}};
\node at (5.1,4.1) {\small{$E_2$}};
\node at (2.9,2.8) {\small{$G_1$}};
\node at (4.9,3.3) {\small{$H_1$}};
\node at (3.8,2.3) {\small{$F$}};
\node at (1.9,1.1) {\small{$G_2$}};
\node at (6,1.5) {\small{$H_2$}};


%
%
%
%

\end{tikzpicture}

Figure 4

\end{center}
\end{minipage}
\end{center}

Figure 3 is a recoloring of Figure 2, and Figure 4 gives an intuitive representation of $A$ and $B$ intersecting, where $V$ is given by the entire rectangle. This is not a proper Venn diagram in the set-theoretic sense, though certain intersections are represented properly, namely $A^{\omega} \cap B$, $B^{\omega} \cap A$ and $A^{\omega} \cap B^{\omega} = I$. 

It remains now only to check that the property i) is fulfilled, \text{i.e.} that the dimensions of the nine summands in our decomposition are uniquely determined by the canonical invariants associated to the pair $(A,B)$. Since any lagrangian subspace of a symplectic subspace has half the dimension of the space within which it is lagrangian, it suffices to show for example that the dimensions of the  subspaces $I$, $E$, $F$, $G_1$ and $H_1$ are uniquely determined.

First, 
$$\dim I = \dim (A^{\omega} \cap B^{\omega}) = k_1$$ 
and the relationships  
\begin{align*}
\dim K & = \dim (A^{\omega} \cap B) + \dim (B^{\omega} \cap A) - \dim I \\
& =  \dim (A^{\omega} \cap B) + [\dim (A^{\omega} \cap B)+ \dim B^{\omega} - \dim A^{\omega}] - \dim I \\
& = 2k_5 + k_3 - k_2 - k_1
\end{align*}
and 
\begin{align*}
\dim W & =  \dim (A^{\omega} + B^{\omega}) \\
& = \dim A^{\omega} + \dim B^{\omega} - \dim  (A^{\omega} \cap B^{\omega}) \\
& = k_2 + k_3 - k_1
\end{align*}
show that $\dim K$ and $\dim W$ are determined.  

Because $E \simeq W/K$ and $F \simeq W^{\omega}/K$ we have 
$$
\dim E = \dim W - \dim K = 2k_2 - 2k_5 
$$
and 
$$
\dim F = \dim (A \cap B) - \dim K = 2k_1 - 2k_3 + 2k_4 - 2k_5
$$

Lastly, $G_1 \simeq (A^{\omega} \cap B) / I$ and $H_1 \simeq (B^{\omega} \cap A) / I$, so $$\dim G_1 = \dim (A^{\omega} \cap B) - \dim I = k_5 - k_1$$ 
and 
$$\dim H_1 = \dim (B^{\omega} \cap A) - \dim I = - k_1 - k_2 + k_3 + k_5$$ 
which proves the property i) and concludes the proof. \qed

\subsection{Elementary types and normal forms}\label{Elementary types and normal forms}

The key to Proposition \ref{main prop} was the decomposition (\ref{inner decomp}), satisfying the properties i) and ii). One may rephrase the construction as follows. We found an $\omega$-orthogonal decomposition 
$$V = V_1 \oplus V_2 \oplus V_3 \oplus V_4 \oplus V_5$$ 
into five symplectic subspaces, such that 

a) the dimensions of these subspaces are uniquely determined by the canonical invariants associated to the coisotropic pair $(A,B)$, and 

b) $A$ and $B$ decompose into direct sums
$$ A = A_1 \oplus A_2 \oplus A_3 \oplus A_4 \oplus A_5 $$
$$ B = B_1 \oplus B_2 \oplus B_3 \oplus B_4 \oplus B_5 $$
such that $A_i \subset V_i$ and $B_i \subset V_i$ for $i = 1,...,5$. 

In other words, we can set $V_1 = D$, $V_2 = E$, $V_3 = F$, etc., and relabel the decompositions
$$A = I \oplus E_1 \oplus G_1 \oplus F \oplus H_1 \oplus H_2$$ 
$$ B = I \oplus E_2  \oplus H_1 \oplus F  \oplus G_1 \oplus G_2$$ 
by setting as $A_i$ as the sum of those summands which lie in $V_i$, i.e. $A_1 = I$, $A_2 = E_1$, $A_3 = F$, $A_4 = G_1$, $A_5 = H_1 \oplus H_2$, and analogously so for $B$.  

Note that for each $i \in \{1,...,5\}$ the subspaces $A_i$ and $B_i$ form a coisotropic pair in $V_i$ of a particularly simple form, each member of the pair being either the entire subspace $V_i$ or a lagrangian subspace therein. Indeed, $A_1 = B_1 = I$ are the same lagrangian subspace of $V_1$, $A_2 = E_1$ and $B_2 = E_2$ form a lagrangian pair in $V_2$, $A_3 = B_3 = F = V_3$, $A_4$ is a lagrangian subspace of $B_4 = G =V_4$, and finally $A_5 = H = V_5$ and $B_5 = H_1$ is lagrangian in this space. We introduce notation for these particularly simple cases of coisotropic pairs. 

\begin{dfn}\label{basic types} A coisotropic pair $(A,B)$ in a symplectic space $V$ is of \textbf{elementary type} if it is one of the following types: 
\begin{description}
\item[$\lambda$:] $A$ and $B$ are lagrangian subspaces, and $A = B$
\item[$\delta$:] $A$ and $B$ are lagrangian subspaces, and $A \cap B = 0$
\item[$\sigma$:] $A = B = V$, i.e. $A$ and $B$ are symplectic
\item[$\mu_B$:] $B = V$ and $A$ is a lagrangian subspace
\item[$\mu_A$:] $A = V$ and $B$ is a lagrangian subspace
\end{description}
We will consider these types ordered as listed and also call them $\tau_1$ through $\tau_5$. 
\end{dfn}

The cases when a coisotropic subspace $C \subset V$ is the entire space or is lagrangian are the two extreme cases of a coisotropic subspace in the sense that they correspond respectively to when $C^{\omega} = 0$ or when $C^{\omega}$ is as large as possible, i.e. $C^{\omega} = C$. The basic types listed above cover all the scenarios when two coisotropics $A$ and $B$ are given by either of these two extremes, except for the possible scenario when $A$ and $B$ are two non-identical lagrangians with non-zero intersection. This case, though, can be split into a ``direct sum'' of the cases $\delta$ and $\lambda$, i.e. it is not ``elementary'' as a type of coisotropic pair. To see this, assume that $A$ and $B$ are such, and let $\tilde A$ and $\tilde B$ be complements of $A \cap B$ in $A$ and $B$ respectively (in particular $\tilde A \cap \tilde B = 0$). Set $W = A + B$ and note that $W^{\omega} = A^{\omega} \cap B^{\omega} = A \cap B \subset W$ because $A$ and $B$ are lagrangian. The subspace $\tilde V := \tilde A \oplus \tilde B$ is such that $\tilde V \oplus (A \cap B) = W$, hence by Lemma \ref{Witt-Artin decomposition} it is symplectic and 
$$ V= \tilde V  \oplus \tilde V ^{\omega}$$
with $A \cap B$ as a lagrangian subspace of $\tilde V ^{\omega}$. With respect to this decomposition of $V$, the coisotropics $A$ and $B$ decompose as $A = \tilde A \oplus A \cap B$ and $B = \tilde B \oplus A \cap B$, where $\tilde A$ and $\tilde B$ are a lagrangian pair in $\tilde V$, i.e. a coisotropic pair of type $\delta$, whereas $A \cap B$, seen as the component of both $A$ and $B$ in $\tilde V^{\omega}$, represents a coisotropic pair in $\tilde V^{\omega}$ of the type $\lambda$.

In the following we make more precise the sense in which a coisotropic pair is the direct sum of smaller coisotropic pairs and in which way the elementary types defined above are indeed elementary.  

\begin{dfn}
Given an $\omega$-orthogonal decomposition of $V$ into a finite number $m \in \mathbb{N}$ of symplectic subspaces
$$ V = \bigoplus_{i=1}^m V_i $$
and given subspaces $A_i$,$B_i \subset V_i$ forming a coiostropic pair in $V_i$ for each $i \in \{1,...,m\}$, we say that $(A,B)$ is the direct sum of the coisotropic pairs $(A_i, B_i)$ if
$$ A = \bigoplus_{i=1}^m A_i \quad \quad \text{and} \quad \quad B = \bigoplus_{i=1}^m B_i $$
Such a direct sum decomposition will be denoted 
$$ (A,B) = \bigoplus_{i}^m (A_i,B_i) $$
\end{dfn}

\begin{dfn} A coisotropic pair $(A,B)$ in $V$ is called \textbf{elementary} if there exists no such direct sum decomposition of $(A,B)$ except as a direct sum of coisotropic pairs of only one of the elementary types $\lambda$, $\delta$, $\sigma$, $\mu_B$ or $\mu_A$. 
\end{dfn}

\begin{prop}\label{elementary types are elementary}
The elementary types $\lambda$, $\delta$, $\sigma$, $\mu_B$ and $\mu_A$ are elementary according to the above definition. 
\end{prop}

\pf Assume that $(A,B)$ is a coisotropic pair of some elementary type 
$$\tau \in \{\lambda, \delta, \sigma, \mu_B, \mu_A\}$$
and that 
$$ (A,B) = \bigoplus_{i}^m (A_i,B_i) $$
is a direct sum decomposition into coisotropics, subordinate to an $\omega$-orthogonal decomposition $V  = \bigoplus V_i$ into symplectic subspaces, i.e. such that $A_i \subset V_i$ and $B_i \subset V_i$ for each $i$. We need to show that each coistropic pair $(A_i,B_i)$ in $V_i$ is of type $\tau$. Because $\tau$ is an elementary type, $A$ is either equal to $V$ or is lagrangian in $V$. If $A = V$, then $A_i = V_i \ \forall i$ for dimension reasons. If $A$ is lagrangian, it is in particular isotropic, and hence each $A_i$ is isotropic in $V$ because $A_i \subset A$. Because $A_i \subset V_i$,  we have $\omega = \omega_{V_i}$ on $A_i$, so $A_i$ is also isotropic in $V_i$. Since $A_i$ is assumed coisotropic in $V_i$, it follows that $A_i$ is lagrangian in $V_i$.  By the same arguments, if $B = V$ then $B_i = V_i \ \forall i$, or if $B$ is lagrangian in $V$ then $B_i$ is lagrangian in $V_i \ \forall i$. It is now clear that if $\tau = \sigma$, then $A_i = B_i = V_i$ for all $i$, so the summand pairs $(A_i,B_i)$ are all also of type $\sigma$. If $\tau =\delta$, then all the $A_i$ and $B_i$ are lagrangian subspaces in their respective $V_i$, and $A \cap B = 0$ implies that $A_i \cap B_i = 0$ for all $i$, so each pair $(A_i,B_i)$ is also of type $\delta$. If $\tau = \lambda$, then similarly the $A_i$ and $B_i$ are lagrangian in $V_i$. To see that here $A_i = B_i \ \forall i$, consider $v \in A = B \subset V$, which has a unique decomposition $v = v_1 + ... + v_m$ with $v_i \in V_i$ for each $i$. Because $v\in A$ and $v \in B$, $v$ also has such unique decompositions with respect to $A = \bigoplus A_i$ and $B = \bigoplus B_i$, but because $A_i,B_i \subset V_i$ for each $i$, these decompositions must coincide with the above decomposition. Hence $v_i \in A_i \cap B_i$ for each $i$. In particular $A = B \subset \bigoplus (A_i \cap B_i)$, which, for dimension reasons, implies $A_i = B_i$ for all $i$.  So each pair $(A_i,B_i)$ is indeed of type $\lambda$ when $(A,B)$ is. Now assume $\tau = \mu_B$. For each $i$, $A_i$ is lagrangian in $V_i$ and $B_i = V_i$, so $(A_i,B_i)$ is also of type $\mu_B$. The case for $\mu_A$ is the same, but with the roles of $A$ and $B$ reversed. \qed

\begin{cor} If a coisotropic pair $(A,B)$ has a direct sum decomposition 
$$ (A,B) = \bigoplus_{i}^m (A_i,B_i) $$
where every coisotropic pair $(A_i,B_i)$ is of the same elementary type, then $(A,B)$ is elementary and of that type. 
\end{cor}

\pf
It suffices to show that $(A,B)$ is of the same type as its summands, since by Proposition \ref{elementary types are elementary} it is then elementary. If the elementary type of the summands is such that the $A_i$ are all lagrangian subspaces of the $V_i$, then the $A_i$ are isotropic subspaces of $V$ and hence their $\omega$-orthogonal sum $A = \bigoplus A_i$ will also be isotropic. Since each $A_i$ has half the dimension of $V_i$, $A$ will have half the dimension of $V$, i.e. it is lagrangian. If on the other hand the elementary type in question is such that  $A_i = V_i$ for each $i$, then clearly $A = V$. The same arguments apply to $B$. Thus the coisotropic pair $(A,B)$ is such that $A$ and $B$ are each either lagrangian or all of $V$ in the same way that their summands $A_i$ and $B_i$ are. It remains only to be sure that when $A$ and $B$ are both lagrangian, they are either identical or such that $A\cap B = 0$, according to whether $A_i = B_i \ \forall i$ or $A_i \cap B_i = 0  \ \forall i$. If $A_i = B_i \ \forall i$ then clearly $A = B$. Assume $A_i \cap B_i = 0 \ \forall i$ and let $v \in A \cap B$. We have a unique decomposition $v = v_1 + ... + v_m$ with $v_i \in V_i$ for all $i$, and because $A = \bigoplus A_i$ and $B = \bigoplus B_i$  are direct sum decompositions subordinate to $V = \bigoplus V_i$, each $v_i$ lies in $A_i \cap B_i = 0$. Hence $v = 0$, and we conclude that $A \cap B = 0$ when $A_i \cap B_i = 0 \ \forall i$.\qed

Proposition \ref{elementary types are elementary} guarantees that the five elementary types of coisotropic pairs are independent of one another in the sense that one cannot express any one of them as a sum of the others. The proof of Proposition \ref{main prop} showed that these basic types are also ``generating'' in the sense that any coisotropic pair decomposes into a direct sum of such elementary types. The corollary implies that one can simplify any direct sum decomposition of a coisotropic pair so that it has only five summands, these summands being of one each of the elementary types. We will call any such five part decomposition an \textbf{{elementary decomposition}}. The following shows that elementary decompositions give a set of invariants for a coisotropic pair $(A,B)$ which are equivalent to the original invariants we associated to such a pair. 

\begin{prop}\label{invariants are equivalent} Let $(A,B)$ be a coisotropic pair in $V$ and let
$$ (A,B) = \bigoplus_{i}^5 (A_i,B_i) $$
be an elementary decomposition subordinate to an $\omega$-orthogonal decomposition $$V= \bigoplus_{i=1}^5 V_i$$ 
ordered such that $(A_i,B_i)$ is of type $\tau_i \in \{ \lambda, \delta, \sigma, \mu_B, \mu_A \}$. Set $n_i := \frac{1}{2} \dim V_i$. Then the 5-tuple 
$$\textbf{n} := (n_1,..., n_5)$$ 
gives a set of invariants (call them \textbf{elementary invariants}) which are equivalent to the canonical invariants 
$$\textbf{k} := (\dim (A^{\omega} \cap B^{\omega}, \ \dim A^{\omega}, \ \dim B^{\omega}, \ \tiny{\frac{1}{2}} \dim V, \  \dim (A^{\omega} \cap B) )$$

\end{prop}

\pf Consider $\textbf{n} = (n_1,...,n_5)$ as a coordinate in the space $\mathcal{N} := \mathbb{Z}_{\scriptscriptstyle{ \geq 0}}^5$ of all possible 5-tuples of elementary invariants (each $V_i$ is symplectic, hence of even dimension), and let $\mathcal{K}$ denote the space of all possible sets of canonical invariants $\textbf{k} = (k_1,...,k_5)$. 

Fix a coisotropic pair $(A,B)$ and fix also an elementary decomposition of this pair, with $A = A_1 \oplus ... \oplus A_5$ and $B = B_1 \oplus ... \oplus B_5$. This gives a 5-tuple $\textbf{n}$. From this $\textbf{n}$ we can obtain the canonical invariants $\textbf{k}$ associated to $(A,B)$ as follows. 

Clearly one has
$$ k_4 =  \tiny{\frac{1}{2}} \dim V = n_1 + n_2 + n_3 + n_4 + n_5 $$
For the remaining invariants, we claim that
$$k_1 = \dim (A^{\omega} \cap B^{\omega}) = n_1 $$
$$ k_2 = \dim A^{\omega} = n_1 + n_2 +  n_4  $$
$$ k_3 = \dim B^{\omega} = n_1 +  n_2 +  n_5 $$
and
$$k_5 = \dim (A^{\omega} \cap B) = n_1 + n_4 $$
To see this, we show 
$$A^{\omega} = A_1 \oplus A_2 \oplus A_4, \quad \quad B^{\omega} =  B_1 \oplus B_2 \oplus B_5$$ 
$$A^{\omega} \cap B = A_1 \oplus A_4, \quad \text{and}  \quad A^{\omega} \cap B^{\omega} = A_1$$
which gives the above formulae for $k_1, k_2,k_3$ and $k_5$ directly. 

For any $a \in A$ we have the decomposition $a = a_1 + ... + a_5$ with $a_i \in A_i$, and for $\tilde a$ also in $A$
\begin{equation}\label{orthogonality splitting}
\omega (a,\tilde a) = \omega_{V_1}(a_1,\tilde a_1) + ... + \omega_{V_5}(a_5,\tilde a_5) = \omega_{V_3}(a_3,\tilde a_3) + \omega_{V_5}(a_5,\tilde a_5)
\end{equation}
because $A_1, A_2$ and $A_4$ are lagrangian in their respective $V_i$. If $\tilde a$ is in $A^{\omega}$, then choosing $a$ as any element  in $A_3$ we find $0 = \omega_{V_3}(a, \tilde a_3)$ and hence $\tilde a_3 \in A_3^{\omega_{V_3}} = 0$, since $A_3$ is symplectic in $V_3$. Similarly one finds $\tilde a_5 = 0$, so $\tilde a \in A_1 \oplus A_2 \oplus A_4$, which shows $A^{\omega} \subset A_1 \oplus A_2 \oplus A_4$. The opposite inclusion $A^{\omega} \supset A_1 \oplus A_2 \oplus A_4$ follows from (\ref{orthogonality splitting}) as well, since for $\tilde a \in A_1 \oplus A_2 \oplus A_4$ and any $a \in A$ we find $\omega (a,\tilde a) = 0$. Arguing analogously one also shows $B^{\omega} = B_1 \oplus B_2 \oplus B_5$.

For the equalities  $A^{\omega} \cap B = A_1 \oplus A_4$ and $A^{\omega} \cap B^{\omega} = A_1$ we use the fact that if $v$ is in $A^{\omega} \cap B$ or $A^{\omega} \cap B^{\omega}$, then in particular $v$ is in $A \cap B$ and hence has a unique decomposition $v = v_1 + ... + v_5$ with $v_i \in A_i \cap B_i \ \forall i$. 

If $v \in A^{\omega} \cap B$, then $v \in A^{\omega} = A_1 \oplus A_2 \oplus A_4$ implies $v_3 = v_5 = 0$. Also, $v_2 \in A_2 \cap B_2 = 0$. Thus $v \in A_1 \oplus A_4$ and $A^{\omega} \cap B \subset A_1 \oplus A_4$ holds. On the other hand, because $A_1 = B_1$ and $A_4 \subset B_4 = V_4$, we have $A_1 \oplus A_4 \subset A^{\omega} \cap B$. 

If $v \in A^{\omega} \cap B^{\omega}$, then not only are $v_3, v_5$ and $v_3$ zero because $A^{\omega} \cap B^{\omega} \subset A^{\omega} \cap B$, but also $v_4 = 0$, because $B^{\omega} = B_1 \oplus B_2 \oplus B_5$ does not contain non-zero summands in $B_4$. Thus $A^{\omega} \cap B \subset A_1$. The opposite in inclusion holds since $A_1 = B_1$ is a summand in the decompositions of both $A^{\omega}$ and $B^{\omega}$.

The equations above describing the $k_i$ in terms of the $n_i$ define a linear map $M: \mathcal{N} \longrightarrow \mathcal{K}$, representable by matrix multiplication with the matrix
$$
M =
\left( \begin{array}{ccccc}
1 & 0 & 0 & 0 & 0 \\
1 & 1 & 0 & 1 & 0 \\
1 & 1 & 0 & 0 & 1 \\
1 & 1 & 1 & 1 & 1 \\
1 & 0 & 0 & 1 & 0 
\end{array} \right)
$$
which is non-singular ($\det M = 1$). Hence $M$ defines an injective map, which means in particular that the numbers $\textbf{n} = (\frac{1}{2} \dim V_1,...,\frac{1}{2} \dim V_5)$ which we associate to an elementary decomposition of a coisotropic pair $(A,B)$ do not depend on the particular elementary decomposition but only depend on the pair $(A,B)$. In other words, $\textbf{n}$ does in fact define a set of invariants for $(A,B)$. The map $M$ is also surjective. Any $\textbf{k} \in \mathcal{K}$ is, by definition, realizable by some coiostropic pair $(A,B)$ and by the proof of Propostion \ref{main prop} this pair has an elementary decomposition; by the above, the invariants $\textbf{n}$ associated to this decomposition are mapped under $M$ to $\textbf{k}$. \qed

To compute the elementary invariants from the canonical invariants one can simply use the inverse of the mapping $M : \textbf{n} \mapsto \textbf{k}$,  
$$
M^{-1} =
\left( \begin{array}{ccccc}
1 & 0 & 0 & 0 & 0 \\
0 & 1 & 0 & 0 & -1 \\
1 & 0 & -1 & 1 & -1 \\
-1 & 0 & 0 & 0 & 1 \\
-1 & -1 & 1 & 0 & 1 
\end{array} \right)
$$
which gives the linear equations for the $n_i$ in terms of the $k_i$:
\begin{align}
n_1 & = k_1\label{n_1 equation} \\
n_2 & = k_2 - k_5 \\
n_3 & = k_1 - k_3 + k_4 - k_5 \\
n_4 & = -k_1 + k_5 \\
n_5 & = -k_1 - k_2 + k_3 + k_5 \label{n_5 equation}
\end{align}
Note that we already nearly explicitly computed these equations in the proof of Proposition \ref{main prop}. 

\begin{cor}\label{ineq only constraint}
The canonical invariants $(k_1,...,k_5)$ are subject only to the five inequalities 
\begin{equation*}
0 \leq k_1 \leq k_5 \leq k_2 \quad \quad k_1 + k_2 \leq k_3 + k_5  \leq k_1 +  k_4 
\end{equation*} 
\end{cor}

\pf
That the $k_i$ must satisfy these inequalities follows from the linear equations (\ref{n_1 equation}) though (\ref{n_5 equation}) for the $n_i$ in terms of the $k_i$ and the fact that $n_i \geq 0 \ \forall i$. The equation for $n_1$ implies $0 \leq k_1$, the equation for $n_2$ gives $k_5 \leq k_2$, the one for $n_3$ gives $k_3 + k_5 \leq k_1 + k_4$, and the inequalities $k_1 \leq k_5$ and $k_1 + k_2 \leq k_3 + k_5$ follow from the equations for $n_4$ and $n_5$. 

To see that these inequalities are the only constraints on the $k_i$, let $\textbf{k} = (k_1,...,k_5)$ be an arbitrary 5-tuple of integers subject only to the above inequalities. We need to show that $\textbf{k}$ is in  $\mathcal{K}$, the set of canonical invariants realizable by a coisotropic pair, which is the image of $M$. In other words we must find a 5-tuple of non-negative integers $\textbf{n} = (n_1,...,n_5)$ such that $M \cdot \textbf{n} = \textbf{k}$, i.e. which solve the linear equations
\begin{align*}
k_1 & = n_1\\
k_2 & = n_1 + n_2 + n_4 \\
k_3 & = n_1 + n_2 + n_5  \\
k_4 & = n_1 + n_2 + n_3 + n_4 + n_5 \\
k_5 & = n_1 + n_4 
\end{align*}
For $k_1 \geq 0$ we choose $n_1 = k_1$ and for $k_5 \geq k_1$ we can always choose $n_4 \geq 0$ such that $k_5 = k_1 + n_4 = n_1 + n_4$. Next, because $k_2 \geq k_5 = n_1 + n_4$, we can choose $n_2 \geq 0$ such that $k_4 = k_5 + n_2 = n_1 + n_2 + n_4$. Thus far $n_1, n_2$ and $n_4$ are fixed and the equations for $k_1, k_2$ and $k_5$ solved. For $k_3$ we have $k_3 \geq k_1 + k_2 - k_5 = n_1 + n_2$, so $n_5$ can be chosen such that $k_3 = n_1 + n_2 + n_5$. Finally, for $k_4 \geq k_3 + k_5 - k_1 = n_1 + n_2 + n_4 + n_5$, an integer $n_3 \geq 0$ is still free to be chosen such that $k_4 = k_3 + k_5 - k_1 + n_3 = n_1 + n_2 + n_3 + n_4 + n_5$ as desired. 
\qed

Using the elementary invariants one can easily construct a normal form $(A_0, B_0)$ for a coisotropic pair $(A,B)$, $\text{i.e.}$ a standardized representative of the equivalence class of $(A,B)$. Let $\textbf{n} = \{n_1,...,n_5\}$ be the elementary invariants of $(A,B)$. We choose $\mathbb{R}^{2n_1} \oplus ... \oplus \mathbb{R}^{2n_5}$ as our model space, equip each summand with the standard symplectic form $\Omega_i$ represented by the $2n_i \times 2n_i$ matrix
$$
\left(  
\begin{array}{cc}
0 & \mathds{1} \\
-\mathds{1} & 0
\end{array}
\right) $$
and give the whole space the direct sum symplectic form $\Omega_1 \oplus ... \oplus \Omega_5$.  Let 
$$(q_1^{i},...,q_{n_i}^{i},p_1^{i},...,p_{n_i}^{i})$$ 
denote the standard coordinates on $\mathbb{R}^{2n_i}$ and denote 
$$Q^{n_i} = \text{span} \{q_1^{i},...,q_{n_i}^{i} \} \text{ and } \ P^{n_i} = \text{span} \{ p_1^{i},...,p_{n_i}^{i} \} $$
Then 
$$ A_0 := Q^{n_1} \oplus Q^{n_2} \oplus \mathbb{R}^{2n_3} \oplus Q^{n_4} \oplus \mathbb{R}^{2n_5}$$
$$ B_0 := Q^{n_1} \oplus P^{n_2} \oplus \mathbb{R}^{2n_3} \oplus \mathbb{R}^{2n_4} \oplus Q^{n_5}$$
defines a normal form for $(A,B)$. By construction $(A_0,B_0)$ is a coisotropic pair such that the elementary invariants of $(A_0,B_0)$ and $(A,B)$ match. Indeed the very definition of $(A_0,B_0)$ gives an elementary decomposition with appropriate dimensions: $(Q^{n_1},Q^{n_1})$ is a coisotropic pair of elementary type $\lambda$ in $\mathbb{R}^{2n_1}$, $(Q^{n_2},P^{n_2})$ a pair of type $\delta$ in $\mathbb{R}^{2n_2}$, and so on. From Proposition \ref{invariants are equivalent} we know that the canonical invariants of $(A,B)$ and $(A_0,B_0)$ match because their elementary invariants do, and by Proposition \ref{main prop} this means that $(A,B) \sim (A_0,B_0)$.

\section{Classification of Linear Canonical Relations}\label{Classification of Linear Canonical Relations}

\subsection{Reduced classification problem}

Let $L \subset V \oplus V^-$ and $\hat L \subset \hat V \oplus \hat V^-$ denote once again canonical relations, and $(A,B)$ and $(\hat A, \hat B)$ the coisotropic pairs giving their respective ranges and domains. We have seen that   for $L$ and $\hat L$ to be equivalent it is necessary that $(A,B)$ and $(\hat A, \hat B)$ be equivalent as coisostropic pairs. Assume now that this is the case, and let $S:V \rightarrow \hat V$ be a symplectic map such that $S(A) = \hat A$ and $S(B) = \hat B$. Furthermore let 
$$
(A,B) = \bigoplus_{i}^5 (A_i,B_i) \quad \text{and }\quad (\hat A,\hat B) = \bigoplus_{i}^5 (\hat A_i,\hat B_i)
$$
be elementary decompositions such that $S(A_i) = \hat A_i$ and  $S(B_i) = \hat B_i$ for each $i$. This situation can be assumed without loss of generality, since the image under $S$ of the decomposition of $(A,B)$ defines an elementary decomposition of $(\hat A, \hat B)$, and this decomposition can in turn can be mapped into any other elementary decomposition of $(\hat A, \hat B)$ by a symplectic map which respects the decompositions. 

The elementary decomposition of $(A,B)$ (and similarly so for $(\hat A, \hat B)$) is subordinate to a decomposition of $V$ of the form
$$ V = (I \oplus J ) \oplus (E_1 \oplus E_1) \oplus F \oplus (G_1 \oplus G_2) \oplus (H_1 \oplus H_2) $$
i.e. $A$ and $B$ decompose as 
$$A = I \oplus E_1 \oplus G_1 \oplus F \oplus H = A_1 \oplus A_2 \oplus A_3 \oplus A_4 \oplus A_5$$ 
$$ B = I \oplus E_2  \oplus H_1 \oplus F  \oplus G = B_1 \oplus B_2 \oplus B_3 \oplus B_4 \oplus B_5$$ 
and we will use the lettered names and indexed names of subspaces interchangeably. Set $A_0 = A_3 \oplus A_5 = F \oplus H$ and $B_0 = B_3 \oplus B_4 = F \oplus G$, \text{i.e.} $A_0$ and $B_0$ are symplectic subspaces such that $A = A^{\omega} \oplus A_0$ and $B = B^{\omega} \oplus B_0$. Recall that $\dim A_0 = \dim B_0$ and $\dim G = \dim H$. From Proposition \ref{Tulczyjew} we know that $L$ and $\hat L$ induce symplectic maps $\phi_L : F \oplus H \rightarrow F \oplus G$ and $\phi_{\hat L} : \hat A_0 \rightarrow \hat B_0$, and that by Proposition \ref{equivalence by parts}, $L \sim \hat L$ if and only if 
\begin{equation}\label{mixed conjugacy condition}
\phi_{\hat L} \circ S \vert_{F \oplus H} = S \vert_{F \oplus G} \circ \phi_L \quad 
\end{equation}
This is equivalent to asking that the diagram 
$$\begin{CD}
F\oplus H @>S\vert_{F\oplus H}>> \hat F \oplus \hat H \\
@V\phi_LVV @VV\phi_{\hat L}V \\
F\oplus G @>S\vert_{F \oplus G}>> \hat F \oplus \hat G
\end{CD}$$
commute. Thus when $L$ and $\hat L$ have matching coisotropic pair invariants, the classification problem reduces to the question of when a map $S$ satisfying this condition exists. 

We consider first two special cases. If $\dim F = 0$, then (\ref{mixed conjugacy condition}) is equivalent to commutativity of the simpler diagram
$$\begin{CD}
H @>S\vert_H>> \hat H \\
@V\phi_LVV @VV\phi_{\hat L}V \\
G @>S\vert_G>> \hat G
\end{CD}$$
Because the invariants of $(A,B)$ and $(\hat A,\hat B)$ match, we can construct an equivalence $S$ (in the sense of coisotropic pairs) by freely choosing symplectic maps $S_1$ through $S_5$ between the symplectic pieces of the decomposition
$$ V = D \oplus E \oplus F \oplus G \oplus H $$
and their corresponding counterparts in the decomposition of $\hat V$, and then taking $S$ as the direct sum of these maps. In particular, in the present case $F = \hat F = 0$ and for any choice of symplectic map $S_5 : H \rightarrow \hat H$ one may choose $S_4 : G \rightarrow \hat G$ as the symplectic map 
$$\phi_{\hat L} \circ S_5 \circ \phi_L^{-1} : G \rightarrow \hat G$$ 
resulting in a map $S$ which both respects the decompositions in $V$ and $\hat V$ and satisfies the condition (\ref{mixed conjugacy condition}). 

The second special case is when $\dim G = \dim H = 0$. Here the condition (\ref{mixed conjugacy condition}) amounts to the condition that $\phi_L$ and $\phi_{\hat L}$ are conjugate via a symplectic map, \text{i.e.} there exists a symplectic map $S_F : F \rightarrow \hat F$ such that the diagram
$$\begin{CD}
F @>S_F>> \hat F \\
@V\phi_LVV @VV\phi_{\hat L}V \\
F @>S_F>> \hat F
\end{CD}$$
commutes. In other words, $L$ and $\hat L$ are equivalent here when $\phi_L$ and $\phi_{\hat L}$, seen as canonical relations in $F \oplus F^-$ and $\hat F  \oplus \hat F^-$ respectively, are equivalent. Thus the classification here is reduced to the symplectic map case. 

For the remainder of this section we now consider the remaining case, \text{i.e.} we assume that $\dim F$, $\dim G$ (and $\dim H$) are all non-zero. This case is at the present moment yet unresolved. We reformulate the problem in coordinates to show what the problem looks like in terms of matrices. Let $2k$ denote the dimension of $F$ and $2l$ the dimension of $G$ and $H$. Let $f = \{f_1,...f_{2k}\}$ be a basis of $F$ and $\Phi_f : \mathbb{R}^{2k} \rightarrow F$ be the corresponding coordinate chart which maps the $i$-th canonical basis vector to $f_i$. Similarly let $g$ and $h$ denote bases of $G$ and $H$, with corresponding charts $\Phi_g$ and $\Phi_h$. The condition (\ref{mixed conjugacy condition}) is equivalent to asking that the diagram
%
\begin{equation*}
\begin{tikzcd}
\mathbb{R}^{2k} \oplus \mathbb{R}^{2l} \arrow{rrr}{\text{$[S_F] \oplus [S_H]$}} \arrow{ddd}[swap]{\text{$[\phi_L]$}} \arrow{dr}{\text{$\Phi_f \oplus \Phi_h$}}
	& 
		& 
			& \mathbb{R}^{2k} \oplus \mathbb{R}^{2l} \arrow{ddd}{\text{$[\phi_{\hat L}]$}} \arrow{dl}[swap]{\text{$\Phi_{\hat f} \oplus \Phi_{\hat h}$}} \\
	& F\oplus H  \arrow{r}{S_F \oplus S_H} \arrow{d}[swap]{\phi_L}
		& \hat F \oplus \hat H \arrow{d}{\phi_{\hat L}} 	
			&  \\
	& F\oplus G  \arrow{r}[swap]{S_F \oplus S_G} 
		& \hat F \oplus \hat G 
			&  \\
\mathbb{R}^{2k} \oplus \mathbb{R}^{2l} \arrow{rrr}[swap]{\text{$[S_F] \oplus [S_G]$}} \arrow{ur}[swap]{\text{$\Phi_f \oplus \Phi_g$}} 
	& 
		& 
			& \mathbb{R}^{2k} \oplus \mathbb{R}^{2l} \arrow{ul}{\text{$\Phi_{\hat f} \oplus \Phi_{\hat g}$}}
\end{tikzcd}
\end{equation*}
commute (the brackets surrounding the maps on the outer rectangle are used to denote the coordinate matrices with respect to the given bases). In other words, finding the maps $S_F$, $S_H$ and $S_H$ as required by (\ref{mixed conjugacy condition}) is equivalent to finding bases $f$,$h$ and $g$ and block diagonal matrices $P \oplus Q = [S_F] \oplus [S_H]$ and $P \oplus R = [S_F] \oplus [S_G]$ such that $P$,$Q$, and $R$ are symplectic and 
\begin{equation}\label{matrix commutation}
[\phi_{\hat L}] = (P \oplus R) \circ [\phi_L] \circ (P \oplus Q)^{-1} 
\end{equation}
If one writes $[\phi_L]$ and $[\phi_{\hat L}]$ as a block matrices 
$$
M =
\left( \begin{array}{ccccc}
M_1 & M_2 \\
M_3 & M_4 
\end{array} \right)
\quad \quad
\hat M =
\left( \begin{array}{ccccc}
\hat M_1 & \hat M_2 \\
\hat M_3 & \hat M_4 
\end{array} \right)
$$
then the commutation condition (\ref{matrix commutation}) reads as 
\begin{equation}\label{matrix equations}
\left( \begin{array}{ccccc}
\hat M_1 & \hat M_2 \\
\hat M_3 & \hat M_4 
\end{array} \right)
=
\left( \begin{array}{ccccc}
PM_1P^{-1} & PM_2Q^{-1} \\
RM_3P^{-1} & RM_4Q^{-1} 
\end{array} \right)
\end{equation}
which gives four matrix equations. A complication here is the fact that, although $M$ represents a symplectic map, the blocks $M_1,...,M_4$ themselves do not (and similarly for $\hat M$). For the problem of finding the equivalence classes of matrices up to the relation (\ref{matrix equations}), the following are possible strategies: 
\begin{enumerate}
\item[i)] put $M$ first into a normal form for symplectic matrices, and then apply the condition (\ref{matrix equations}) 
\item[ii)] find a normal form for the condition (\ref{matrix equations}) without any symplectic assumptions, and then apply symplectic constraints as a second step
\end{enumerate}
In the case of either i) or ii), an apparent issue is the following. The normal form given by Gutt \cite{Gutt}, for example, arises from a decomposition of a symplectic map into a direct sum of symplectic maps on generalized eigenspaces. This decomposition need not ``respect" the splitting $A_0 = F \oplus H$, and thus one is faced again with the problem of how one of the normal form blocks ``intertwines" the spaces $F$, $H$ and $G$. Furthermore it is a priori unclear which of the normal form blocks leave $F$ invariant and map $H$ to $G$ (and hence are unaffected by the condition (\ref{matrix equations})) and which do not. Because the literature on normal forms for symplectic matrices is diverse, it is possible that a different normal form than the one given by Gutt in \cite{Gutt} would be more amenable to the condition (\ref{matrix equations}). The study of other normal form structures would thus be one natural step in the further study of this problem. 

Besides the approach of considering the problem in coordinate matrices, one might obtain a simplification of the of the classification associated to (\ref{mixed conjugacy condition}) by first considering invariants of this subproblem which are built from the dimensions and ranks of the subspaces given by the intersections of $\phi_L(F)$ and $\phi_L(H)$ with $F$ and $G$. Additionally, a computation of the decomposition of $L$ as a direct sum of Towber's basic types of linear relations might give insight into the structure of the induced relation $\phi_L$. These steps, yet incomplete, are at present omitted from this exposition. 

\subsection{Normal forms}\label{Normal forms}

We present one possible way of representing a linear canonical relation which reflects its decomposition into linear canonical relations of simpler types, and is also amenable to the yet to be completed full classification. As a subspace in $V \oplus V^-$, any linear canonical relation $L$ is specified by a basis for this subspace, which in split coordinates can be written as
$$
\left( \begin{array}{c}
 v_1  \\
 w_1 	
\end{array} \right)
,... \ ,
\left( \begin{array}{c}
 v_n  \\
 w_n 	
\end{array} \right)
$$
with $v_i, w_i \in V$ and $v_i L w_i$ for each $i \in \{1,...,n\}$. These pairs can be arranged vertically as the columns of a matrix having $2n$ rows and $n$ columns, ordered to reflect the decompositions of the coisotropic subspaces which are the domain and range of $L$.  Clearly, for any element of $v \in \text{ker}(L) = A^{\omega}$ one has $vL0$, and for $v \in \text{hal}(L)=B^{\omega}$ holds $0Lv$. Because $\text{dom}(L) = A$ and $\text{ran}(L) = B$ decompose as 
$$
A = A^{\omega} \oplus A_0 = (I \oplus E_1 \oplus G_1) \oplus (F \oplus H)
$$
$$
B = B^{\omega} \oplus B_0 = (I \oplus E_2 \oplus H_1) \oplus (F \oplus G)
$$
we can represent $L$ in matrix form as 
$$
\left( \begin{array}{cccccccc}
 I & E_1	 &  G_1 	& F 	 		& H		& 	0 	& 0		&  0 \\
0 & 0		 &  0   	& \phi_L(F) 	& \phi_L(H) & I	& E_2 	& H_1 
\end{array} \right)
$$
where each letter stands in for a basis of the subspace it denotes. Recall that in the decomposition 
$$
V = D \oplus E \oplus F \oplus G \oplus H
$$
$I \subset D$ is lagrangian subspace, and $(E_1,E_2)$ form a transverse lagrangian pair in $E$. One sees that  the block columns corresponding to $I$ in the representation of $L$ above represent a linear canonical relation $L_D$ in $D \oplus D^-$ given by the block matrix 
$$
\left( \begin{array}{cc}
 I & 0	  \\
0 & I	
\end{array} \right)
$$
and similarly the columns corresponding to $E_1$ and $E_2$ represent a linear canonical relation $L_E$
$$
\left( \begin{array}{cc}
 E_1 & 0	  \\
0 & E_2
\end{array} \right)
$$
in $E \oplus E^-$. Setting $V_0 = F \oplus G \oplus H$ we thus have a decomposition
$$
V \oplus V^- = (D \oplus D^-) \oplus (E \oplus E^-) \oplus (V_0 \oplus V_0^-)
$$
and 
$$
L = L_D \oplus L_E \oplus L_0
$$ 
where $L_0 \subset V_0 \oplus V_0^-$ is the linear canonical relation with ker$(L) = G_1$,  hal$(L) =  H_1$ and with induced symplectic map $\phi_{L_0} : F \oplus H \rightarrow F \oplus G$ which coincides with $\phi_L$. In terms of the matrix representation of $L$ we write
$$
\left( \begin{array}{cc}
 I & 0	  \\
0 & I	
\end{array} \right)
\oplus
\left( \begin{array}{cc}
 E_1 & 0	  \\
0 & E_2
\end{array} \right)
\oplus
\left( \begin{array}{cccc}
 G_1 &F 	 		& H	& 	0 \\
 0 & \phi_{L_0}(F) 	& \phi_{L_0}(H) & H_1	
\end{array} \right)
$$
Finally, using the notation from the classification and normal forms for coisotropic pairs one can choose as a canonical normal form the ``block matrix" given by
$$
\left( \begin{array}{cc}
Q^{n_1}& 0	  \\
0 	& Q^{n_1} 
\end{array} \right)
\oplus
\left( \begin{array}{cc}
Q^{n_2} & 0	  \\
0 & P^{n_2}
\end{array} \right)
\oplus
\left( \begin{array}{cccc}
Q^{n_4} & \mathbb{R}^{2(n_3 + n_4)} & 	0 \\
 0 &  [\phi_{L_0}] & Q^{n_5}	
\end{array} \right)
$$
where one is implicitly assuming the use of the canonical bases in the spaces $Q^{n_1}$, $P^{n_2}$, \text{etc.}, and $[\phi_{L_0}]$ denotes here a yet to be determined general normal form for $\phi_{L_0}$.

\end{document}